\documentclass[12pt, letterpaper]{elsarticle}
\usepackage{amsmath}
\usepackage{amssymb}
\usepackage{graphicx}
\usepackage{multirow}
\usepackage{latexsym}
\usepackage{epic}
\usepackage{url}
\usepackage{soul}
\usepackage{afterpage}

\usepackage{bbm}
\usepackage{setspace}
\usepackage[norelsize,ruled, lined,linesnumbered]{algorithm2e}
\usepackage{enumitem}
\usepackage{cleveref}
\usepackage[toc]{appendix}

\usepackage{tikz}
\usepackage{tikz,fullpage}
\usepackage{multirow}
\usepackage{pgf}
\usetikzlibrary{arrows,automata}
\usepackage{tkz-berge}
\usepackage[position=top]{subfig}
\usepackage{amsthm}

\usepackage[left=0.9in,top=0.9in,right=0.9in,bottom=0.9in,nohead]{geometry}

\makeatletter
\newtheorem*{rep@theorem}{\rep@title}
\newcommand{\newreptheorem}[2]{%
\newenvironment{rep#1}[1]{%
 \def\rep@title{#2 \ref{##1}}%
 \begin{rep@theorem}}%
 {\end{rep@theorem}}}
\makeatother

\newreptheorem{example}{Example}
\newtheorem{theorem}{Theorem}

\newtheorem{lemma}{Lemma}

\newtheorem{remark}{Remark}
\newtheorem{example}{Example}

\usepackage{xcolor}
\usepackage[colorinlistoftodos,prependcaption,textsize=tiny]{todonotes}
\usepackage{xargs}

\usepackage[symbol]{footmisc}

\begin{document}
\onehalfspace
\begin{frontmatter}
\title{An approach to the distributionally robust shortest path problem \footnote{© 2021. This manuscript version is made available under the CC-BY-NC-ND 4.0 license http://creativecommons.org/licenses/by-nc-nd/4.0/}}

\author[label1]{Sergey S.~Ketkov}
\author[label2]{Oleg A.~Prokopyev\footnote[2]{Corresponding author. Email: droleg@pitt.edu; phone: +1-412-624-9833.}}
\author[label1]{Evgenii P.~Burashnikov}

\address[label1]{Laboratory of Algorithms and Technologies for Networks Analysis, National Research University \\Higher School of Economics (HSE), Bolshaya Pecherskaya st., 25/12, Nizhny Novgorod, 603155, Russia}
\address[label2]{Industrial Engineering, University of Pittsburgh, 1048 Benedum Hall, Pittsburgh,  PA 15261, USA}

\begin{abstract}
In this study we consider the shortest path problem, where the arc costs are subject to distributional uncertainty. Basically, the decision-maker attempts to minimize her worst-case expected loss over an ambiguity set (or a family) of candidate distributions that are consistent with the decision-maker's initial information. The ambiguity set is formed by all distributions that satisfy prescribed linear first-order moment constraints with respect to subsets of arcs and individual probability constraints with respect to particular arcs. 
 Under some additional assumptions the resulting distributionally robust shortest path problem (DRSPP) admits equivalent robust and mixed-integer programming (MIP) reformulations. The robust reformulation is shown to be $NP$-hard, whereas the problem without the first-order moment constraints is proved to be polynomially solvable. We perform numerical experiments to illustrate the advantages of the considered approach; we also demonstrate that the MIP reformulation of DRSPP can be solved effectively using off-the-shelf solvers.
\end{abstract}

\begin{keyword}
shortest path problem; distributionally robust optimization; mixed-integer programming; polyhedral un\-cer\-tain\-ty
\end{keyword}

\end{frontmatter}
\onehalfspace
\section{Introduction} \label{sec: intro}
The \textit{shortest path problem} (SPP) has been attracting much interest both theoretically and computationally since the early 1950s \cite{Dijkstra1959, Ahuja1988}. Being one of the classical network optimization problems it finds various applications in transportation, planning, network interdiction and design; see, e.g., \cite{Ford1958, Johnson1978, Israeli2002, Irnich2005} and the references therein.

\looseness-1Consider a weighted directed connected graph $G = (N, A, \mathbf{c})$, where $N$ and $A$ denote its sets of nodes and directed arcs, respectively. With each arc $a \in A$ we associate a nonnegative cost $c_a$, that is, $\mathbf{c} = \{c_a: a \in A\}$. We refer to $s \in N$ and $t \in N$ as the source and destination nodes, respectively. Recall that the standard deterministic problem of finding an $s-t$ path of the minimum total cost is known to be polynomially solvable, e.g., by dynamic programming algorithms of Dijkstra and Bellman-Ford \cite{Dijkstra1959, Bellman1958}. We also refer to \cite{Ahuja1988, Bazaraa2011} for a more comprehensive discussion.

However, in practice the decision-maker often does not know the nominal arc costs/travel times in advance. In fact, uncertain factors such as variability of travel times, path capacity variation may significantly influence the quality of routing decisions; see, e.g., \cite{Bertsimas2003, Mudchanatongsuk2008}. The modeling approach for data uncertainty depends on a concrete application, but in general consists of the following two major principles.

On the one hand, a \textit{robust optimization approach} represents unknown costs/travel times\footnote{We use the terms ``costs'' and ``travel times'' interchangeably.} through uncertainty sets, i.e., the cost vector $\mathbf{c}$ is assumed to belong to some uncertainty~set~$\mathcal{S}$. Then a particular measure of robustness is optimized across all possible realizations of costs $\mathbf{c}~\in~\mathcal{S}$; see surveys \cite{Bental1998, Bental2002, Bental2009, Bertsimas2011}. Despite a great modeling power robust solutions assume no distributional knowledge and thus, potentially provide overly conservative decisions.

On the other hand, a \textit{stochastic programming approach} assumes that the cost vector $\mathbf{c}$ is governed by some known probability distribution $\mathbb{Q}^0$; see, e.g., \cite{Ruszczynski2003, Shapiro2009}, which is referred to as the \textit{nominal distribution}. In this case one may optimize some risk measure under the specified distribution $\mathbb{Q}^0$. Nevertheless, it is often argued, see, e.g., \cite{Wiesemann2014}, that fitting a single candidate distribution to the available information potentially leads to biased optimization results with poor out-of-sample performance. What is probably more important, the distribution of the cost vector is often not known to the decision-maker in advance; see \cite{Zhang2017} and the references therein.

Alternatively, a \textit{distributionally robust optimization approach} represents the uncertainty by an \textit{ambiguity set} (or a family) $\mathcal{Q}$ of probability distributions that are compatible with the decision-maker's initial information; see, e.g., related studies in \cite{Wiesemann2014, Delage2010, Goh2010, Bertsimas2018}. The idea is to optimize some utility function across the constructed family of probability distributions, i.e., with respect to $\mathbb{Q} \in \mathcal{Q}$. This approach attempts to balance between the lack of distributional information and the complete knowledge of the underlying distribution. In particular, Wiesemann~et~al.~\cite{Wiesemann2014} introduce a unified approach to solving distributionally robust convex optimization problems. In this paper we adopt the optimization techniques proposed in \cite{Wiesemann2014} to the \textit{shortest path problem with distributional uncertainty}.

\subsection{Related literature} \label{subsec: lit review}
This section provides an overview of the most related literature and is organized as follows. First, we briefly outline robust optimization and stochastic programming models for the shortest path problem. Next, we describe existing formulations of SPP under distributional uncertainty.

The literature on the robust shortest path problem (RSPP) is vast, see, e.g., a survey in \cite{Kasperski2016}. Various types of uncertainty sets including polyhedral, discrete or budgeted ones, are used to model a variability of the cost vector; see \cite{Dokka2017}. In particular, the robust shortest path problem with interval data is considered in \cite{Bertsimas2003, Montemanni2004}.

Typically, robust optimization methods focus on the control of a \textit{conservatism level} of the proposed solutions.
For example, Bertsimas~et~al.~\cite{Bertsimas2003} introduce a parameter $\Gamma$, which can be used to limit the maximal number of components of vector $\mathbf{c}$ that deviate from their nominal values. Hence, by varying $\Gamma$ the decision-maker is able to control her level of protection against uncertainty in a more sophisticated manner. The robust optimization approach proposed in \cite{Bertsimas2003} preserves polynomial solvability, while most of the robust versions of the shortest path problem are $NP$-hard in general; see, e.g., \cite{Dokka2017, Yu1998, Zielinski2004}. Naturally, the robust optimization approach may lead to suboptimal decisions when some distributional information is available to the decision-maker; this observation is also validated numerically in our computational experiments, see the discussion in Section~\ref{sec: comp study}.

\looseness-1Next, we refer to \cite{Hall1983, Bertsekas1991, Nie2009, Chen2010, Pan2013} for stochastic programming models related to SPP. Typically such models assume that the nominal distribution of the cost vector is known to the decision-maker and seek optimal paths with respect to some predefined reliability criterion. Such problems are usually computationally challenging since even evaluation of the objective function for a fixed decision requires calculation of a multidimensional integral \cite{Shapiro2008}. Thus, some approximation methods such as sample average approximation and discretization have been developed; see, e.g., \cite{Pan2013, Verweij2003}. At the same time, in practice we frequently encounter a lack of data to reconstruct the nominal distribution of arc costs/travel times; we refer to \cite{Ng2011} for more details.

An interesting version of the stochastic shortest path problem is studied by Rinehart~et~al.~\cite{Rinehart2011}. In their setting side information about arc costs is incorporated through a specialized conditional expectation model. The authors provide analytic information bounds for the performance of their model and highlight some practical insights when the travel times are Gaussian. 

Also, several studies consider the \textit{distributionally robust shortest path problem} (DRSSP). For example, one line of research focuses on distributionally robust mixed-integer programming problems; we refer to \cite{Bertsimas2004, Bertsimas2006, Natarajan2011} and the review by Li~et~al.~\cite{Li2014}. The major contributions of these research studies can be summarized as follows. (\textit{i}) The authors consider max-min problems, where decisions are implemented \textit{after} realization of uncertainty; also, (\textit{ii}) the ambiguity sets include marginal moment, cross-moment sets as well as some generalizations (that account for both moment and cross-moment information). The authors exploit conic duality to obtain finite-dimensional reformulations of the corresponding distributionally robust problems. (\textit{iii}) It is shown in \cite{Bertsimas2006} that for the marginal moment ambiguity set (the set with a predefined marginal mean and variance for each component of the cost vector) the resulting distributionally robust formulation is polynomially solvable, as long as the original MIP problem is also polynomially solvable. Furthermore, for the cross-moment ambiguity set (which accounts for the mean and covariance of the cost vector) the problem is proved to be $NP$-hard even for a class of linear programs \cite{Natarajan2011}. (\textit{iv}) The authors define persistency of a decision variable, i.e., the probability that it is contained in the optimal solution. It is shown in \cite{Bertsimas2006} that under the marginal moment model persistency can be computed efficiently as a byproduct of the basic solution procedure.

Additionally, a number of studies consider min-max formulations of DRSPP where routing decisions are made \textit{here-and-now}, prior to the realization of uncertainty. Thus, Gavriel~et~al.~\cite{Gavriel2012} explore dynamic programming-based approximations for the risk-averse shortest path problem with distributional uncertainty. The arc travel times are assumed to be independent random variables with a given mean and variance.

Zhang~et~al.~\cite{Zhang2017} and Cheng~et~al.~\cite{Cheng2016} propose min-max versions of DRSPP, which consider correlation between arc costs. In their settings the distributional ambiguity set accounts for information about the support, mean and covariance matrix of uncertain parameters. Guided by the work in \cite{Delage2010}, they view DRSPP as a mixed-integer semidefinite programming problem. Since the resulting problem is computationally difficult, a sequence of semidefinite relaxations is considered and the tightness of the obtained bounds is demonstrated numerically.



\looseness-1Recently, Jaillet~et~al.~\cite{Jaillet2016} have proposed an index-based approach to a class of distributionally robust routing problems. The authors introduce a Requirements Violation (RV) index, which quantifies the risk associated with an uncertain parameter violating specified upper and lower limits. The uncertain parameter is given by individual arc travel times/uncertain demands or their linear combination. Specifically, RV index resolves a trade-off between two extreme situations, i.e., a situation where the limits are met almost surely (resulting in support constraints) and a situation where the limits are satisfied only in expectation (resulting in the worst-case expectation constraints). Then the problem of minimizing cumulative RV index is formulated and the corresponding solution techniques based on Benders decomposition are developed.

The travel times/uncertain demands in \cite{Jaillet2016} are assumed to be either independent or inherently dependent via a linear regression model. Importantly, the proposed ambiguity sets are moment-based, but yield a confidence region for the mean (and possibly higher-order moments) of each uncertain parameter.

\looseness-1Finally, some data-driven approaches for modeling the distributional uncertainty are also available in the literature. In such models the distribution of the uncertain parameters is only observable through a finite training dataset. Then a family of distributions is designed based on a distance metric in the space of probability distributions ``centered'' at the empirical distribution of the training samples. We refer to \cite{Wang2020} for a model of DRSPP based on Wasserstein balls and to \cite{Zhu2018} for a formulation of network sensor allocation problem based on phi-divergences.

\subsection{The approach and contribution in this study} \label{subsec: contribution}
\looseness-1The primary goal of this paper is to develop an approach to capture distributional uncertainty in the context of the shortest path problem. In the view of discussion above, we assume that the cost vector $\mathbf{c}$ is subject to distributional uncertainty. We explore a min-max version of DRSPP, i.e., optimize \textit{here-and-now} decisions implemented before the realization of uncertainty.

The key feature of our model is that instead of moment-based ambiguity sets, which account for the first- and second-order moments explicitly, we use standardized ambiguity sets proposed by Wiesemann~et~al.~\cite{Wiesemann2014}.
This prerequisite allows us to construct the ambiguity sets from partially observable and unreliable data. In addition, our distributional constraints can be strengthened with a new information available and therefore, can be adopted in sequential decision-making settings; see our further discussion below and in Section \ref{subsec: modeling the ambiguity set}.

Next, we briefly discuss our construction of the distributional constraints. Furthermore, we provide some practical insights and motivation behind the choice of both the objective function and ambiguity set in our modeling approach.

Let $\mathbb{Q}$ be a joint distribution of the random vector $\mathbf{c}$. Then for a given vector $\mathbf{b} \in \mathbb{R}^{D_0}$ and real-valued matrix $\mathbf{B} \in \mathbb{R}^{D_0 \times |A|}$ we introduce linear expectation constraints of the form:
\begin{equation} \label{eq: linear expectation constraints}
\mathbb{E}_{\mathbb{Q}} \{\mathbf{B} \mathbf{c}\} \leq \mathbf{b},
\end{equation}
where $D_0 \in \mathbb{Z}_+$ denotes the number of expectation constraints.

By leveraging~(\ref{eq: linear expectation constraints}) the decision-maker can limit the cumulative expected cost of any $s-t$ path or any subset of arcs $A' \subseteq A$. This idea finds applications in online learning frameworks, where the feedback received by the decision-maker is typically limited; see, e.g., \cite{Dani2008, Bubeck2012}. For example, she may have access to some historical data containing observations of a \textit{total cost} with respect to particular routes in the network. In Section \ref{subsec: modeling the ambiguity set} we discuss how this information can be exploited to construct linear expectation constraints of the form (\ref{eq: linear expectation constraints}) subject to a prescribed confidence level.
In particular, if we introduce auxiliary random variables $\mathbf{u} \in \mathbb{R}^{D_0}$ such that $\mathbf{u} := \mathbf{b} - \mathbf{B} \mathbf{c}$, then linear expectation constraints (\ref{eq: linear expectation constraints}) coincide with those considered by Wiesemann~et~al.~\cite{Wiesemann2014}.

Let $\mathbb{Q}_a \in \mathcal{Q}_0(\mathbb{R})$, $a \in A$, be a marginal distribution induced by $\mathbb{Q}$. Then guided by \cite{Wiesemann2014} for each particular arc $a \in A$ we introduce \textit{individual support} and \textit{quantile constraints}:
\allowdisplaybreaks\begin{gather}
\mathbb{Q}_a\{c_a \in [l_a, u_a]\} = 1 \label{eq: support constraint} \\
\mathbb{Q}_a\{c_a \in [l^{(i)}_a, u^{(i)}_a]\} \in [\underline{q}^{(i)}_a, \overline{q}^{(i)}_a], \mbox{ } i \in \mathcal{D}_a,\label{eq: quantile constraints}
\end{gather}
where $0 \leq l_a \leq u_a < \infty$; $[l^{(i)}_a, u^{(i)}_a] \subseteq [l_a, u_a], \mbox{ } i \in \mathcal{D}_a := \{1, \ldots, D_a\}$, is a set of $D_a \in \mathbb{Z}_+$ subintervals; $\underline{q}^{(i)}_a$ and $\overline{q}^{(i)}_a$ specify the probability that random cost~$c_a$ belongs to the $i$-th subinterval, where $0 \leq \underline{q}^{(i)}_a \leq \overline{q}^{(i)}_a \leq 1$.
\begin{remark} \label{remark support}\upshape
Note that support constraints (\ref{eq: support constraint}) can be represented in the form of (\ref{eq: quantile constraints}). Hence, without loss of generality, we assume that for each $a \in A$ and $i = D_a$ the constraint

$$\mathbb{Q}_a\{c_a \in [l^{(i)}_a, u^{(i)}_a]\} \in [\underline{q}^{(i)}_a, \overline{q}^{(i)}_a]$$
\noindent is a support constraint with $l^{(i)}_a = l_a$, $u^{(i)}_a = u_a$ and $\underline{q}^{(i)}_a = \overline{q}^{(i)}_a = 1$.
\vspace{-10mm}\flushright$\square$
\end{remark}

Next, we briefly discuss the intuition behind quantile constraints in the form of (\ref{eq: quantile constraints}). On the one hand, these constraints can be constructed from interval-censored or unreliable data. In other words, to model (\ref{eq: quantile constraints}) it is necessary to verify whether the arc cost $c_a$ belongs to a specified subinterval, i.e., $c_a \in [l^{(i)}_a, u^{(i)}_a]$ for some $i \in \mathcal{D}_a$. Meanwhile, pointwise observations of the arc costs can be either unreliable or unavailable to the decision-maker; we refer the reader to Section~\ref{subsec: modeling the ambiguity set} for more details.

On the other hand, by leveraging quantile constraints (\ref{eq: quantile constraints}) the decision-maker can substantially improve the quality of a distributionally robust optimization model. The role of individual probability constraints is examined, e.g., in \cite{Clare2012, Chen2017}, with applications to air traffic flow management. As a remark, it is rather straightforward to verify that distributional constraints (\ref{eq: quantile constraints}) can provide a high quality approximation of the marginal distribution $\mathbb{Q}_a$, $a \in A$. For example, one can imagine a situation where the subintervals are sufficiently small, mutually disjoint and cover the support interval $[l_a, u_a]$, whereas the upper and lower quantiles, namely, $\underline{q}^{(i)}_a$ and $\overline{q}^{(i)}_a$, $i \in \mathcal{D}_a$, coincide.

For each particular arc $a \in A$ we consider separate probability constraints. In fact, probability constraints with respect to subsets of arcs must satisfy so called \textit{nesting condition} for the distributionally robust counterpart (i.e., the second-level optimization problem in the min-max formulation) to be a tractable convex problem; see Theorems 1 and 2 in \cite{Wiesemann2014}. This condition is quite restrictive since it requires for the confidence sets (in our case, the subintervals $[l^{(i)}_a, u^{(i)}_a]$, $i \in \mathcal{D}_a$, $a \in A$) to be either disjoint or nested within each other. However, given a specialized structure of the objective function and quantile constraints in our problem formulation we avoid using the nesting condition, but preserve tractability of the second-level problem.

With respect to the objective function, our study focuses on minimizing of the \textit{worst-case expected loss} incurred by the decision-maker.
From the modeling perspective our choice of the objective function can be justified as follows.


First, it is outlined in \cite{Miller2000} that in some applications, especially those of a repetitive nature, it may be sufficient to find the paths with minimal expected travel time. In our setting repetitive decisions arise naturally, if the decision-maker learns the nominal distribution by trial and errors through multiple decision epochs. Specifically, she may refine some distributional information by implementing her solutions sequentially several times; see Section \ref{subsec: modeling the ambiguity set} for modeling the distributional constraints. Besides, the paths with the least expected cost are often used in intelligent transportation and in-vehicle route guidance systems; see, e.g., \cite{Fu1998}.

Under the expected loss criterion we identify a form of the worst-case distribution and provide an efficient solution procedure for DRSPP without linear expectation constraints. With respect to the general case, we obtain a linear mixed-integer programming (MIP) reformulation of DRSPP in contrast to the semidefinite MIP formulations proposed in \cite{Zhang2017, Cheng2016}. In our computational experiments we demonstrate that the proposed MIP formulation can be solved rather effectively using off-the-shelf solvers.



Note that linear MIP reformulations of combinatorial optimization problems with objective uncertainty can also be derived, if one uses optimized certainty equivalent risk measures combined with specified expectation constraints; see, e.g., \cite{Hanasusanto2016}. We briefly discuss the applicability of some other objective criteria to our problem setting in Section~\ref{sec: comp study}.   
 
As a remark, bounds on the covariance matrix used in \cite{Delage2010, Cheng2016} exploit the mean of the cost vector $\mathbf{c}$. However, if the expected costs are known, then the objective function value in our setting is fully defined.
Therefore, we ``handle'' correlation information between travel times by leveraging linear expectation constraints (\ref{eq: linear expectation constraints}). 
 
In summary, this paper describes a rather simple approach for solving the shortest path problem under distributional uncertainty. To the best of our knowledge, the vast majority of ambiguity sets proposed in the related literature (including moment-based sets, Wasserstein balls and etc.)~are highly sensitive to both incompleteness and unreliability of historical data. For instance, an accurate estimation of the correlation matrix (as it is proposed in~\cite{Zhang2017}) typically requires a sufficiently large number of pointwise observations with respect to each component of vector~$\mathbf{c}$; see, e.g., \cite{Bickel2008}. Contrariwise, we demonstrate that the distributional constraints considered in this study can be constructed based on interval-censored observations for particular arcs and observations of the total cost for some subsets of arcs.

In the view of discussion above the contributions and the remaining structure of the paper can be summarized as follows:
\begin{itemize}
\item In Section \ref{sec: problem} we formulate DRSPP and discuss how to obtain the distributional constraints (\ref{eq: linear expectation constraints}) and (\ref{eq: quantile constraints}) from real-data observations.
\item\looseness-1 In Section \ref{sec: DR SP without EC} we prove that DRSPP without expectation constraints can be solved in polynomial time by retrieving optimal costs for each particular arc and solving a deterministic shortest path problem. Also, we provide some structural observations on the form of the worst-case distribution.
\item Section \ref{sec: general} provides a robust formulation of DRSPP with a polyhedral uncertainty set. The resulting bilevel optimization problem turns out to be $NP$-hard in general and we describe its single-level MIP reformulation.
\item In Section \ref{sec: comp study} we conduct a numerical study using the proposed approach and evaluate its out-of-sample performance. Additionally, we provide a comparison with standard robust and distributionally robust optimization techniques. 
\item Section~\ref{sec: conclusion} concludes the paper and outlines possible directions for future research.
\end{itemize}

\noindent \textbf{Notation.} All vectors and matrices are labelled by bold letters. Arc $a \in A$ adjacent to nodes $v_1, v_2 \in N$ is denoted as $(v_1, v_2)$. Let $\mathcal{P}_{st}(G)$ be the set of all simple directed paths from $s$ to $t$ in the network $G$. Any path $P \in \mathcal{P}_{st}(G)$ is given by a sequence of arcs $(s,v_1),(v_1,v_2),\ldots,(v_{|P|-1},t)$, which we introduce as $\{s \rightarrow v_1 \rightarrow \ldots \rightarrow v_{|P|-1} \rightarrow t\}$ for convenience. For a subset of arcs $A' \subseteq A$ we define a subgraph of $G$ induced by this subset of arcs as $G[A'] := (N, A', \mathbf{c}')$, where, in particular, $\mathbf{c}' := \{c_a, a \in A'\}$.

We use $\mathbbm{1}\{Z\}$ as an indicator of event~$Z$. 
The uniform distribution on an interval $[l,u]$ is referred to as $\mathcal{U}(l,u)$. Finally, denote by $\mathcal{M}_+(\mathbb{R}^{k})$ and $\mathcal{Q}_0(\mathbb{R}^{k})$ the spaces of all nonnegative measures and probability distributions on $\mathbb{R}^{k}$ for some $k \in \mathbb{Z}_+$, respectively.

\section{Problem formulation} \label{sec: problem}
\subsection{Distributionally robust shortest path problem} \label{subsec: DR SP}
As outlined in Section \ref{sec: intro} the cost vector $\mathbf{c}$ is assumed to be a nonnegative random vector governed by some unknown joint distribution $\mathbb{Q} \in \mathcal{Q}_0(\mathbb{R}^{|A|})$. With each arc $a \in A$ we associate a marginal probability distribution $\mathbb{Q}_a \in \mathcal{Q}_0(\mathbb{R})$ induced by $\mathbb{Q}$.

The joint distribution $\mathbb{Q}$ is supposed to belong to an ambiguity set $\mathcal{Q}$ comprised of all probability distributions that satisfy linear expectation constraints (\ref{eq: linear expectation constraints}) and individual quantile constraints (\ref{eq: quantile constraints}). That is,
\begin{equation} \label{ambiguity set}
\begin{gathered}
\mathcal{Q} = \Big\{\mathbb{Q} \in \mathcal{Q}_0(\mathbb{R}^{|A|}) \mbox{: } \mathbb{E}_{\mathbb{Q}} \{\mathbf{B} \mathbf{c}\} \leq \mathbf{b};\\
\qquad \qquad \qquad \qquad \qquad \quad \; \; \mathbb{Q}_a\{c_a \in [l^{(i)}_a, u^{(i)}_a]\} \in [\underline{q}^{(i)}_a, \overline{q}^{(i)}_a] \quad \forall i \in \mathcal{D}_a, \mbox{ } a \in A
\Big\}
\end{gathered}
\end{equation}

For each node $i \in N$ we refer to $RS_i$ ($FS_i$) as the set of the arcs directed out of (into) node~$i$. Denote by $\mathbf{y} \in \{0,1\}^{|A|}$ a path-incidence vector and introduce the standard flow-balance constraints \cite{Ahuja1988} as:
\begin{equation} \label{eq: flow-balance constraints}
\mathbf{y} \in Y = \{\mathbf{y} \in \{0,1\}^{|A|}: \mathbf{G}\mathbf{y} = \mathbf{g}\},
\end{equation}
where $\mathbf{G} \in \{-1, 0, 1\}^{|N|\times|A|}$ and $\mathbf{g} \in \{0,1\}^{|N|}$. Specifically, for each $i \in N$
$$g_i = \begin{cases}
1, \mbox{ if } i = s\\
-1, \mbox{ if } i = f \\
0, \mbox{ otherwise}
\end{cases}$$
Furthermore,
$$G_{ij} = \begin{cases} 1, \mbox{ if } j \in RS_i\\
-1, \mbox{ if } j \in FS_i\\
0, \mbox{ otherwise }\\
\end{cases}$$

Then the \textit{distributionally robust shortest path problem} (DRSPP) is formulated as follows:
\begin{align} \label{DRO formulation} \nonumber \tag{\textbf{F1}}
\min_{\mathbf{y} \in Y} \max_{\mathbb{Q} \in \mathcal{Q}} \mathbb{E}_{\mathbb{Q}}\{\mathbf{c}^\top \mathbf{y}\}
\end{align}

That is, we minimize the worst-case expected loss of the decision-maker, i.e., the worst-case expected path cost, across all probability distributions consistent with the decision-maker's prior information. Henceforth, we need the following modeling assumptions:\\
\textbf{A1.} For each $a \in A$ there exists a marginal distribution $\mathbb{Q}_a \in \mathcal{Q}_0(\mathbb{R})$ such that
$$\mathbb{Q}_a\{l^{(i)}_a \leq c_a \leq u^{(i)}_a\} \in (\underline{q}^{(i)}_a, \overline{q}^{(i)}_a),$$
whenever $\underline{q}^{(i)}_a < \overline{q}^{(i)}_a$, $i \in \mathcal{D}_a$.\\
\textbf{A2.} For each $a \in A$ and any pair of subintervals in (\ref{eq: quantile constraints}), namely, $[l^{(i_1)}_a, u^{(i_1)}_a]$ and $[l^{(i_2)}_a, u^{(i_2)}_a]$, $i_1, i_2 \in \mathcal{D}_a$, we have $l^{(i_1)}_a \neq u^{(i_2)}_a$ and $l^{(i_2)}_a \neq u^{(i_1)}_a$.



Assumption \textbf{A1} guarantees existence of a probability distribution that satisfies quantile constraints (\ref{eq: quantile constraints}) as strict inequalities, if interval $[\underline{q}^{(i)}_a, \overline{q}^{(i)}_a]$ is non-degenerate. Additionally, it allows us to exploit the strong duality results for the moment problems in Section \ref{sec: DR SP without EC}; we refer the reader to \cite{Wiesemann2014, Shapiro2001} for a more comprehensive discussion.

Next, Assumption \textbf{A2} stipulates that the inner optimization problem in (\ref{DRO formulation}) has a finite maximum; see the proof of Lemma \ref{lemma2} in Section \ref{sec: DR SP without EC} and Example \ref{example A2}. We emphasize that Assumption \textbf{A2} is not restrictive since, if, e.g., $l^{(i_1)}_a = u^{(i_2)}_a$ for some $i_1, i_2 \in \mathcal{D}_a$, then a sufficiently small perturbation of the endpoints makes this assumption satisfied.

\subsection{Data-driven approach for modeling the ambiguity set} \label{subsec: modeling the ambiguity set}
Next, we discuss how to construct the family of distributions (\ref{ambiguity set}) from real-data observations. For simplicity we assume that, if the decision-maker traverses through $s-t$ paths multiple times, then the joint distribution $\mathbb{Q}$ of the cost vector is fixed across all decision epochs.

As outlined in Section \ref{subsec: contribution}, linear expectation constraints in the form of (\ref{eq: linear expectation constraints}) model a situation, where the decision-maker observes only the total path cost in each decision epoch. Thus, suppose that a path $P \in \mathcal{P}_{st}(G)$ is traversed by the decision-maker $r \in \mathbb{Z}_+$ times. We refer to $\boldsymbol{\xi}^{(P)} \in \mathbb{R}^r$ as a vector comprised of $r$ i.i.d. observations of the total path cost $\sum_{a \in P} c_a$. Using support constraints (\ref{eq: support constraint}) observe that:

$$l^{(P)} : = \sum_{a \in P} l_a \leq \xi_i^{(P)} \leq u^{(P)} := \sum_{a \in P} u_a, \quad \forall i \in \{1, \ldots, r\}$$
Furthermore, Hoeffding inequality \cite{Boucheron2003} for the sum $\sum_{i = 1}^r \xi_i^{(P)}$ of $r$ bounded i.i.d. random variables implies that for any $\varepsilon > 0$ one has:
\begin{equation} \label{eq: azuma-hoeffding 1}
\mathbb{Q}_P\Big\{|\mathbb{E}_{\mathbb{Q}}\{\sum_{a \in P} c_a\} - \frac{1}{r}\sum_{i = 1}^r \xi_i^{(P)}| \geq \varepsilon\Big\} \leq 2 \exp\Big(\frac{-2r \varepsilon^2}{(u^{(P)} - l^{(P)})^2}\Big),
\end{equation}
where $\mathbb{Q}_P \in \mathcal{Q}_0(\mathbb{R})$ is a distribution of the empirical mean $\frac{1}{r}\sum_{i = 1}^r \xi_i^{(P)}$. Hence, with high probability the following expectation constraints hold:
\begin{equation} \label{eq: to model linear expectation constraints}
\frac{1}{r}\sum_{i = 1}^r \xi_i^{(P)} - \varepsilon \leq \mathbb{E}_{\mathbb{Q}}\{\sum_{a \in P} c_a\} \leq \frac{1}{r}\sum_{i = 1}^r \xi_i^{(P)} + \varepsilon
\end{equation}
Specifically, $\varepsilon$ is defined by setting the right-hand side of (\ref{eq: azuma-hoeffding 1}) equal to a prescribed confidence level. Also, note that instead of complete $s-t$ paths one may consider any nonempty subset of arcs $A' \subseteq A$.

Then support constraints (\ref{eq: support constraint}) can be derived from some physical limitations, i.e., the arc costs/travel times are typically bounded depending upon a concrete application. We also refer to \cite{Delage2010} for a construction of the support constraints based on empirical data.

Finally, we discuss how to construct individual quantile constraints (\ref{eq: quantile constraints}). In the sequel, we fix $a \in A$ and consider the quantile constraints associated with the marginal distribution $\mathbb{Q}_a$. Denote by $\boldsymbol{\xi}^{(a)} \in \mathbb{R}^r$ a vector of $r \in \mathbb{Z}_+$ i.i.d. observations of random cost $c_a$ and pick a subinterval $[l'_a, u'_a] \subseteq [l_a, u_a]$. Furthermore, we define
\begin{equation} \nonumber
\chi_i^{(a)} = \begin{cases}
1, \mbox{ if } \xi_i^{(a)} \in [l'_a, u'_a] \\
0, \mbox{ otherwise,}
\end{cases}
\end{equation}
where $i \in \{1, \ldots, r\}$. Then $\chi_i^{(a)} \in \{0, 1\}$ is a Bernoulli random variable with an unknown probability of success $q_a \in [0, 1]$, that is,
$$q_a = \mathbb{Q}_a\{l'_a \leq c_a \leq u'_a\}$$
Let $\overline{\mathbb{Q}}_a \in \mathcal{Q}(\mathbb{R})$ be a distribution of the empirical mean $\frac{1}{r}\sum_{i = 1}^r \chi_i^{(a)}$. Using Hoeffding inequality for arbitrary $\varepsilon > 0$ observe that:
\begin{equation} \label{eq: azuma-hoeffding 2}
\overline{\mathbb{Q}}_a\Big\{|q_a - \frac{1}{r}\sum_{i = 1}^r \chi_i^{(a)}| \geq \varepsilon\Big\} \leq 2 \exp(-2r \varepsilon^2)
\end{equation}
As a result, with high probability we have:\\
\begin{equation} \label{eq: to model quantiles}
\mathbb{Q}_a\{l'_a \leq c_a \leq u'_a\} = q_a \in [\frac{1}{r}\sum_{i = 1}^r \chi_i^{(a)} - \varepsilon; \frac{1}{r}\sum_{i = 1}^r \chi_i^{(a)} + \varepsilon]
\end{equation}
where parameter $\varepsilon$ depends on a prescribed confidence level and thus, can be defined from (\ref{eq: azuma-hoeffding 2}).

In general, there are no theoretical limitations with respect to the choice of subinterval $[l'_a, u'_a]$. However, one can adopt some practical limitations, e.g., $c_a \leq u'_a$ may indicate that the decision-maker has completed a task in time (and is untimely, otherwise). In Section \ref{sec: comp study} we provide some analysis on how the width of the subintervals influences the optimal solution of DRSPP (\ref{DRO formulation}) in the absence of expectation constraints.


As a result, following the discussion in Section \ref{sec: intro} the construction of linear expectation constraints relies on observations of a total cost with respect to some subsets of arcs. Likewise, the construction of quantile constraints does not require precise knowledge of data observations, but only indicates whether an observation belongs to some predefined subinterval. Therefore, the observations can be collected indirectly, e.g., by checking whether the arc costs/travel times satisfy specified upper and lower limits. We also refer to \cite{Sun2007, Ferson2007} for a comprehensive analysis of interval-censored data.

As a remark, instead of Hoeffding inequality one can employ more advanced measure concentration results from \cite{Boucheron2003}. Furthermore, as outlined in \cite{Wiesemann2014} to guarantee a specified confidence level for the ambiguity set $\mathcal{Q}$ one can adopt confidence levels of the individual constraints by using Bonferroni's inequality \cite{Birge2011}.

\section{Model without expectation constraints} \label{sec: DR SP without EC}
In this section we examine DRSPP (\ref{DRO formulation}) without linear expectation constraints (\ref{eq: linear expectation constraints}). We prove that the resulting problem can be solved in polynomial time. More precisely, it is tackled by solving a particular linear programming problem for each $a \in A$ and a single deterministic shortest path problem.

Hereafter, we suppose that the ambiguity set of probability distributions is given by:
\begin{equation} \label{ambiguity set without expectation constraints}
\begin{gathered}
\widetilde{\mathcal{Q}} := \Big\{\mathbb{Q} \in \mathcal{Q}_0(\mathbb{R}^{|A|}) \mbox{: }
\mathbb{Q}_a\{c_a \in [l^{(i)}_a, u^{(i)}_a]\} \in [\underline{q}^{(i)}_a, \overline{q}^{(i)}_a] \quad \forall i \in \mathcal{D}_a, \mbox{ } a \in A
\Big\}
\end{gathered}
\end{equation}
Consider the following DRSPP without linear expectation constraints:
\begin{equation} \label{DRO formulation without EC} \nonumber \tag{\textbf{F1}$^\prime$}
\min_{\mathbf{y} \in Y} \max_{\mathbb{Q} \in \widetilde{\mathcal{Q}}} \mathbb{E}_{\mathbb{Q}}\{\mathbf{c}^\top \mathbf{y}\}
\end{equation}

First, leveraging the structure of (\ref{ambiguity set without expectation constraints}) we show that optimization problem (\ref{DRO formulation without EC}) can be partitioned into $|A|$ individual moment problems with respect to each particular arc $a \in A$ and then resolved as a deterministic shortest path problem. The following result holds.
\begin{lemma} \label{lemma1}
Let
\begin{equation} \label{marginal ambiguity set}
\widetilde{\mathcal{Q}}_a := \Big\{\mathbb{Q}_a \in \mathcal{Q}_0(\mathbb{R}) \mbox{: }
\mathbb{Q}_a\{c_a \in [l^{(i)}_a, u^{(i)}_a]\} \in [\underline{q}^{(i)}_a, \overline{q}^{(i)}_a] \quad \forall i \in \mathcal{D}_a
\Big\}
\end{equation}
Suppose that $\mathbf{y}^* \in Y$ and $\mathbb{Q}^* \in \mathcal{Q}_0(\mathbb{R}^{|A|})$ is an optimal solution of (\ref{DRO formulation without EC}). In particular, let $\mathbb{Q}_a^* \in \mathcal{Q}_0(\mathbb{R})$, $a \in A$, be the marginal distributions induced by $\mathbb{Q}^*$.
Then
\begin{itemize}
\item for each $a \in A$ the worst-case expected cost $\mathbb{E}_{\mathbb{Q}_a^*}\{c_a\}$ coincides with the optimal objective function value of the following individual moment problem:
\begin{equation} \label{individual primal problem}
\max_{\mathbb{Q}_a \in \widetilde{\mathcal{Q}}_a} \mathbb{E}_{\mathbb{Q}_a} \{c_a\}
\end{equation}
\item an optimal path-incidence vector $\mathbf{y}^*$ can be attained by solving a deterministic shortest path problem of the form:
\begin{equation} \label{deterministic shortest path problem}
\min_{\mathbf{y} \in Y} \sum_{a \in A} \mathbb{E}_{\mathbb{Q}_a^*} \{c_a\} y_a
\end{equation}
\end{itemize}
\begin{proof} \upshape

Note that the set of constraints in (\ref{ambiguity set without expectation constraints}) can be partitioned into $|A|$ non-overlapping subsets, i.e., the quantile constraints for each particular arc $a \in A$. Thus, taking into account the linearity of expectation DRSPP (\ref{DRO formulation without EC}) can be equivalently reformulated as:
\begin{equation} \nonumber
\min_{\mathbf{y} \in Y} \sum_{a \in A} \Big(\max_{\mathbb{Q}_a \in \widetilde{\mathcal{Q}_a}} \mathbb{E}_{\mathbb{Q}_a}\{c_a\}\Big) y_a
\end{equation}
and the result follows.
\end{proof}
\end{lemma}

Next, we apply the duality theory to solve the \textit{individual moment problem} (\ref{individual primal problem}) for each particular arc; see \cite{Shapiro2001}. For simplicity of our further exposition we need the following preprocessing step. For each arc $a \in A$ from the \textit{baseline} set $[l^{(i)}_a, u^{(i)}_a]$, $i \in \mathcal{D}_a$, of subintervals we form a set $[L^{(j)}_a, U^{(j)}_a]$, $j \in \mathcal{W}_a := \{1, \ldots W_a\}$, of $W_a \in \mathbb{Z}_+$ \textit{elementary subintervals} \cite{Berg1997}.

Specifically, consider a list of distinct interval endpoints, that is,
$$\{l^{(1)}_a, u^{(1)}_a, l^{(2)}_a, u^{(2)}_a, \ldots l^{(D_a)}_a= l_a, u^{(D_a)}_a = u_a\}$$
and sort them in a nondecreasing order. Regions of the resulting partitioning of interval $[l_a, u_a]$ are referred to as elementary subintervals and denoted by $[L^{(j)}_a, U^{(j)}_a]$, $j \in \mathcal{W}_a$. For instance, a baseline set of subintervals $$\{[20, 60], [30, 70], [0, 100]\}$$ is split into a set $$\{[0, 20], [20, 30], [30, 60], [60, 70], [70, 100]\}$$ of elementary subintervals.

For any $j \in \mathcal{W}_a$ we denote by $\mathcal{D}_a(j) \subseteq \mathcal{D}_a$ indices of the baseline subintervals contained in the elementary subinterval $[L^{(j)}_a, U^{(j)}_a]$. Analogously, for any $i \in \mathcal{D}_a$ we denote by $\mathcal{W}_a(i) \subseteq \mathcal{W}_a$ indices of the elementary subintervals contained in the baseline subinterval $[l^{(i)}_a, u^{(i)}_a]$.

Finally, we note that the overall complexity of preprocessing is given by $O(D_a \log D_a)$ for each $a \in A$ and the number of nonempty elementary subintervals does not exceed $2D_a - 1$ by construction. Now, we are ready to introduce an equivalent linear programming reformulation of the individual moment problem~(\ref{individual primal problem}).

\begin{lemma} \label{lemma2} Optimization problem (\ref{individual primal problem}) for fixed $a \in A$ can be equivalently reformulated as:

\allowdisplaybreaks\begin{subequations} \label{individual linear programming problem}
\begin{align}
& \max_{\boldsymbol{\delta}_a} \sum_{j \in \mathcal{W}_a} U^{(j)}_a \delta_{aj} \\
\mbox{\upshape s.t. }
& \delta_{aj} \geq 0 \quad \forall j \in \mathcal{W}_a\\
& \sum_{j \in \mathcal{W}_a} \delta_{aj} = 1 \\
& \underline{q}^{(i)}_a \leq \sum_{j \in \mathcal{W}_a(i)} \delta_{aj} \leq \overline{q}^{(i)}_a, \quad \forall i \in \mathcal{D}_a \setminus \{1\}
\end{align}
\end{subequations}

\begin{proof}
Optimization problem (\ref{individual primal problem}) for fixed $a \in A$ coincides with the following moment problem:
\vspace{-5mm}
\begin{subequations} \label{individual primal problem 2}
\begin{align}
& \max_{\mathbb{\mu}}\int_{l_a}^{u_a} c_a \mbox{ d}\mu(c_a) \\
\mbox{\upshape s.t. }
& \mu \in \mathcal{M}_+(\mathbb{R}) \label{cons: individual primal problem 1}\\
& \int_{l_a}^{u_a} \mathbbm{1}\Big\{l^{(i)}_a \leq c_a \leq u^{(i)}_a\Big\} \mbox{ d}\mu(c_a) \leq \overline{q}^{(i)}_a \quad \forall i \in \mathcal{D}_a \label{cons: individual primal problem 2}\\
& \int_{l_a}^{u_a} \mathbbm{1}\Big\{l^{(i)}_a \leq c_a \leq u^{(i)}_a\Big\} \mbox{ d}\mu(c_a) \geq \underline{q}^{(i)}_a \quad \forall i \in \mathcal{D}_a, \label{cons: individual primal problem 3}
\end{align}
\end{subequations}
Assumption \textbf{A1} implies that the strong duality holds; see Proposition 3.4 in \cite{Shapiro2001}. More precisely, individual moment problem (\ref{individual primal problem}) is a particular case of the moment problem considered by Wiesemann~et~al. \cite{Wiesemann2014};   see equation (4) in \cite{Wiesemann2014}.

Denote by $\mathbf{k}_a \in \mathbb{R}_+^{D_a}$ and $\mathbf{h}_a \in \mathbb{R}_+^{D_a}$ dual variables corresponding to the primal constraints (\ref{cons: individual primal problem 2}) and (\ref{cons: individual primal problem 3}), respectively. Then the dual reformulation of (\ref{individual primal problem 2}) is given~by:
\begin{subequations} \allowdisplaybreaks \label{individual robust dual problem}
\begin{align}
& \min_{\mathbf{k}_a,\mathbf{h}_a} \sum_{i \in \mathcal{D}_a}(\overline{q}^{(i)}_a k_{ai} - \underline{q}^{(i)}_a h_{ai}) \\
\mbox{\upshape s.t. }
& k_{ai} \geq 0, \mbox{ } h_{ai} \geq 0 \quad \forall i \in \mathcal{D}_a\\
& \sum_{i \in \mathcal{D}_a} \mathbbm{1}\Big\{l^{(i)}_a \leq c_a \leq u^{(i)}_a\Big\}(k_{ai} - h_{ai}) - c_a \geq 0 \quad \forall c_a \in [l_a, u_a] \label{cons: robust dual}
\end{align}
\end{subequations}
Note that the set of constraints (\ref{cons: robust dual}) is satisfied if and only if
\begin{equation} \allowdisplaybreaks \label{eq: robust counterpart}
\min_{c_a \in [l_a, u_a]} \Big\{\sum_{i \in \mathcal{D}_a} \mathbbm{1}\Big\{l^{(i)}_a \leq c_a \leq u^{(i)}_a\Big\}(k_{ai} - h_{ai}) - c_a\Big\} \geq 0
\end{equation}
Assume that $c_a \in [L^{(j)}_a, U^{(j)}_a]$ for some $j \in \mathcal{W}_a$.
Then using Assumption \textbf{A2} observe that:
$$\sum_{i \in \mathcal{D}_a} \mathbbm{1}\Big\{l^{(i)}_a \leq c_a \leq u^{(i)}_a\Big\}(k_{ai} - h_{ai}) = \sum_{i \in \mathcal{D}_a(j)}(k_{ai} - h_{ai})$$
and the minimum value in the left-hand side of (\ref{eq: robust counterpart}) is achieved at $c_a = U^{(j)}_a$.

Otherwise, if Assumption \textbf{A2} is violated, then there exist $i_1, i_2 \in \mathcal{D}_a$ and $j_0 \in \mathcal{W}_a$ such that $u^{(i_1)}_a = l^{(i_2)}_a = U^{(j_0)}_a$. In this case, the function
\begin{equation} \nonumber
\phi(c_a) := \sum_{i \in \mathcal{D}_a} \mathbbm{1}\Big\{l^{(i)}_a \leq c_a \leq u^{(i)}_a\Big\}(k_{ai} - h_{ai}) - c_a
\end{equation}
is discontinuous at $c_a = U^{(j_0)}_a$, and the minimum value in the left-hand side of (\ref{eq: robust counterpart}) does not necessarily exist.

Taking the minimum across the set of all elementary subintervals, i.e., with respect to $j \in \mathcal{W}_a$, results in the following dual reformulation:
\allowdisplaybreaks\begin{subequations} \label{individual dual problem}
\begin{align}
& \min_{\mathbf{k}_a,\mathbf{h}_a} \sum_{i \in \mathcal{D}_a}(\overline{q}^{(i)}_a k_{ai} - \underline{q}^{(i)}_a h_{ai}) \\
\mbox{\upshape s.t. }
& k_{ai} \geq 0, \mbox{ } h_{ai} \geq 0 \quad \forall i \in \mathcal{D}_a\\
& \min_{j \in \mathcal{W}_a}\Big\{\sum_{i \in \mathcal{D}_a(j)}(k_{ai} - h_{ai}) - U^{(j)}_a\Big\} \geq 0 \label{individual dual problem min}
\end{align}
\end{subequations}

Observe that (\ref{individual dual problem}) can be viewed as a linear programming (LP) problem. Specifically, constraint (\ref{individual dual problem min}) is equivalent to $W_a$ linear constraints of the form:
\begin{equation} \nonumber
\sum_{i \in \mathcal{D}_a(j)}(k_{ai} - h_{ai}) - U^{(j)}_a \geq 0, \quad \forall j \in \mathcal{W}_a
\end{equation}
Recall that for any $a \in A$ we have $\overline{q}^{(1)}_a = \underline{q}^{(1)}_a = 1$ by construction. Then by the standard LP duality we derive a dual formulation of (\ref{individual dual problem}) and the result follows.

\end{proof}
\end{lemma}
\begin{remark} \label{remark strong duality} \upshape
Let $\mathbb{Q}^*_a \in \mathcal{Q}_0(\mathbb{R})$ and $\boldsymbol{\delta}^*_a$ be optimal solutions of optimization problems (\ref{individual primal problem}) and (\ref{individual linear programming problem}), respectively, for some fixed $a \in A$. Then by the strong duality the \textit{worst-case expected cost}, i.e., the expected cost under the worst-case distribution $\mathbb{Q}^*_a$, is given by:

\begin{equation} \label{eq: strong duality!}
\mathbb{E}_{\mathbb{Q}^*_a} \{c_a\} =
\sum_{j \in \mathcal{W}_a}U^{(j)}_a \delta^*_{aj}
\end{equation}

Furthermore, Lemma \ref{lemma2} provides us an intuition behind the construction of the worst-case distribution. Thus, it is a discrete distribution, which takes values at the endpoints of elementary subintervals, $U^{(j)}_a$, $j \in \mathcal{W}_a$.
 \vspace{-9.5mm}\flushright$\square$
\end{remark}

The next example explains some motivation behind Assumption \textbf{A2}. Specifically, we demonstrate that, if Assumption \textbf{A2} is not satisfied, then the maximum in the moment problem (\ref{individual primal problem}) may not be attained.

\begin{example} \label{example A2} \upshape For some $a \in A$ let $l_a = 0$, $u_a = 100$ and $D_a = 3$. Consider two baseline subintervals, i.e., $$[l^{(1)}_a, u^{(1)}_a] = [0, 50] \mbox{ and } [l^{(2)}_a, u^{(2)}_a] = [50, 100]$$ that do not satisfy Assumption \textbf{A2}. Assume that
$$\mathbb{Q}_a\{c_a \in [0, 50]\} = \mathbb{Q}_a\{c_a \in [50, 100]\} = 0.5$$
Then the worst-case distribution is a discrete distribution such that:
\begin{equation} \nonumber
c_{a} = \begin{cases}
50 - \varepsilon, \mbox{ with probability } 0.5 \\
100, \mbox{ with probability } 0.5
\end{cases}
\end{equation}
for arbitrarily small $\varepsilon > 0$. We conclude that the maximum in (\ref{individual primal problem}) does not exist.
\vspace{-9mm}\flushright$\square$
\end{example}

The next theorem states a key theoretical result of this section.

\begin{theorem} \label{theorem1}
Optimization problem (\ref{DRO formulation without EC}) is polynomially solvable.
\begin{proof} \upshape
In fact, Lemmas \ref{lemma1}, \ref{lemma2} and Remark \ref{remark strong duality} imply that optimization problem (\ref{DRO formulation without EC}) can be tackled by solving $|A|$ linear programming problems (\ref{individual linear programming problem}) and a single deterministic shortest path problem (\ref{deterministic shortest path problem}). Observe that each LP has $W_a$ variables and $2D_a + W_a + 1$ constraints. Since the number of elementary subintervals $W_a$ does not exceed $2D_a - 1$ by construction we conclude that optimization problem (\ref{individual linear programming problem}) is polynomially solvable for each~$a \in A$.
\end{proof}
\end{theorem}


In particular, Theorem \ref{theorem1} mimics the results of Bertsimas~et~al.~\cite{Bertsimas2006} for the distributionally robust max-min problems with a marginal moment ambiguity set. Naturally, if the underlying combinatorial optimization problem is polynomially solvable, then its distributionally robust version with quantile constraints is also polynomially solvable. Also, if quantile constraints (\ref{eq: quantile constraints}) coincide with support constraints (\ref{eq: support constraint}), then $c_a \in [l_a, u_a]$, $a \in A$, results in the expectation constraints of the form $\mathbb{E}_{\mathbb{Q}_a}\{c_a\} \in [l_a, u_a]$. Hence, in this particular case (\ref{DRO formulation without EC}) is equivalent to the min-max robust shortest path problem with interval data \cite{Kasperski2016}. The next example demonstrates that exploiting partial distributional information potentially improves the quality of myopic robust solutions.

\begin{figure}
\centering
\begin{tikzpicture}[scale=0.75,transform shape]
\Vertex[x=0,y=0]{1}
\Vertex[x=6,y=2.5]{2}
\Vertex[x=6,y=-2.5]{3}
\Vertex[x=12,y=0]{4}
\tikzstyle{LabelStyle}=[fill=white,sloped]
\tikzstyle{EdgeStyle}=[post]
\tikzstyle{EdgeStyle}=[post, bend left, thick, double = black, double distance = 0.5pt]
\Edge[label=$[0\mbox{,}100]$](1)(2)
\tikzstyle{EdgeStyle}=[post, bend left]
\Edge[label=$[1\mbox{,}101]$](2)(4)
\Edge[label=$[0\mbox{,}100]$](2)(3)
\tikzstyle{EdgeStyle}=[post, bend right]
\Edge[label=$[0\mbox{,}100]$](1)(3)
\Edge[label=$[0\mbox{,}100]$](3)(4)
\end{tikzpicture}
\onehalfspacing
\caption{\footnotesize The network used in Example \ref{ex: Example 1}. The cost range $[l_a, u_a]$ is depicted inside each arc $a \in A$. The arc $a'$ is highlighted in bold.}
\label{fig: Example 1}
\end{figure}
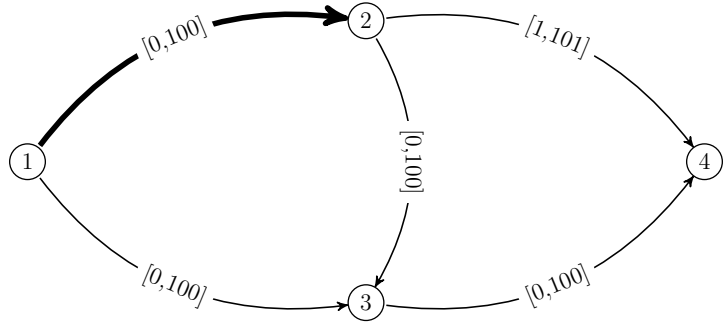

\begin{example} \label{ex: Example 1} \upshape
Consider the network depicted in Figure \ref{fig: Example 1}. Let $s = 1$ and $t = 4$. Suppose that besides the support information for arc $a' = (1,2)$ we know that its cost $c_{a'}$ exceeds $70$ with probability of at most $0.1$. Therefore, we add the following quantile constraint:
$$\mathbb{Q}_{a'}\{70 \leq c_{a'} \leq 100\} \in [0, 0.1]$$

By Lemma \ref{lemma2} and Remark \ref{remark support} the worst-case expected cost of arc $a'$ is obtained by solving the following individual linear programming problem:
\begin{subequations} \label{example individual linear programming problem}
\begin{align}
& \min_{\mathbf{\delta}_{a'}} \quad 70 \delta_{a' 1} + 100 \delta_{a' 2} \\
\mbox{\upshape s.t. }
& \delta_{a' i} \geq 0 \quad \forall i \in \{1, 2\}\\
& \delta_{a' 1} + \delta_{a' 2} = 1 \\
& 0 \leq \delta_{a' 2} \leq 0.1
\end{align}
\end{subequations}

Here, $D_{a'} = W_{a'} = 2$ and the elementary subintervals are given by $[0, 70]$ and $[70, 100]$. Also, $\mathcal{W}_{a'}(1) = \{2\}$ and $\mathcal{W}_{a'}(2) = \{1, 2\}$. The worst-case distribution $\mathbb{Q}^*_{a'}$ with respect to arc~$a'$ is a discrete distribution such that:
\begin{equation} \nonumber
c_{a'} = \begin{cases}
100, \mbox{ with probability } 0.1 \\
70, \mbox{ with probability } 0.9
\end{cases}
\end{equation}
Hence, the worst-case expected cost of $a'$ is given by:
$$\mathbb{E}_{\mathbb{Q}^*_{a'}}\{c_{a'}\} = 0.9 \times 70 + 0.1 \times 100 = 73$$

\begin{table}
\footnotesize
  \centering
  \onehalfspacing
    \begin{tabular}{c c c c c c}
    \hline
    $s-t$ path &Naive robust approach & Distributionally robust approach & Nominal solution \\
    \hline
    $\{1 \rightarrow 2 \rightarrow 4\}$ & 201 & \textbf{174} & \textbf{88.5}\\
    $\{1 \rightarrow 3 \rightarrow 4\}$ & \textbf{200} & 200 & 100\\
    $\{1 \rightarrow 2 \rightarrow 3 \rightarrow 4\}$ & 300 & 273 & 137.5\\
    \hline
    \end{tabular}
  \caption{The nominal expected cost and expected cost under both naive robust and distributionally robust uncertainty for each path $P \in \mathcal{P}_{st}(G)$ in Example \ref{ex: Example 1}. For the former approach all the arc costs are set to their upper bounds.}
  \label{tab: example comparison}
\end{table}

Next, we define a \textit{nominal distribution} $\mathbb{Q}^0 \in \mathcal{Q}_0(\mathbb{R}^{|A|})$ of the cost vector $\mathbf{c}$, which is compatible with the decision-maker's initial information. That is, for $a \neq a'$ suppose that the induced marginal distribution $\mathbb{Q}^0_a \in \mathcal{Q}_0(\mathbb{R})$ is uniform, i.e., $c_a \sim \mathcal{U}(l_a, u_a)$. Additionally, for $a = a'$ assume that:
\begin{equation} \nonumber
c_{a'} \sim \begin{cases}
\mathcal{U}(0, 70), \mbox{ with probability } 0.95 \\
\mathcal{U}(70, 100), \mbox{ with probability } 0.05
\end{cases}
\end{equation}
Therefore, the \textit{nominal expected cost} of $a'$ is given by:
$$\mathbb{E}_{\mathbb{Q}^0_{a'}}\{c_{a'}\} = 0.95 \times 35 + 0.05 \times 85 = 37.5$$

In Table \ref{tab: example comparison} for each $s - t$ path we report its nominal expected cost as well as the expected cost under both robust and distributionally robust uncertainty. The optimal values are in bold. In fact, a naive robust optimization approach sets all the arc costs to their upper bounds. Observe that the optimal distributionally robust solution, in turn, dominates the optimal robust solution in terms of the nominal expected cost. \vspace{-10mm}\flushright$\square$
\end{example}

Example \ref{ex: Example 1} illustrates that robust solutions may suffer from underspecification, especially, when some partial distributional information is available to the decision-maker. In Section~\ref{sec: comp study} we extend the methodology behind Example \ref{ex: Example 1} to a class of randomly generated problem instances with both quantile and linear expectation constraints. The next section provides robust and mixed-integer programming reformulations of DRSPP~(\ref{DRO formulation}).

\section{General case} \label{sec: general}
\subsection{Equivalent robust reformulation} \label{subsec: robust reformulation}
We exploit the results outlined in the previous section to reformulate DRSPP (\ref{DRO formulation}) as an instance of the robust shortest path problem (RSPP) with a polyhedral uncertainty set. Furthermore, we show that the resulting problem is $NP$-hard in general and describe a linear MIP reformulation of (\ref{DRO formulation}).



Suppose that a path-incidence vector $\mathbf{y} \in Y$ is fixed. Now, we analyze the inner optimization problem in (\ref{DRO formulation}), that is,
\begin{equation} \label{lower-level optimization problem}
\max_{\mathbb{Q} \in \mathcal{Q}} \sum_{a \in A} \mathbb{E}_{\mathbb{Q}_a}\{c_a\} y_a
\end{equation}
The key observation is that optimization problem (\ref{lower-level optimization problem}) can be reformulated as a linear programming problem in terms of the marginal expectations $\mathbb{E}_{\mathbb{Q}_a}\{c_a\}$, $a \in A$.
Specifically, we show that for each $a \in A$ there exists a \textit{surjective mapping} from a set of marginal probability distributions
\begin{equation} \label{set of marginal distributions}
\widetilde{\mathcal{Q}}_a = \Big\{\mathbb{Q}_a \in \mathcal{Q}_0(\mathbb{R}): \mathbb{Q}_a\{c_a \in [l^{(i)}_a, u^{(i)}_a]\} \in [\underline{q}^{(i)}_a, \overline{q}^{(i)}_a] \quad \forall i \in \mathcal{D}_a\Big\}
\end{equation}
onto a set formed by linear expectation constraints
\begin{equation} \label{new linear expectation constraints}
c^{min}_a \leq \mathbb{E}_{\mathbb{Q}_a}\{c_a\} \leq c^{max}_a,
\end{equation}
for some $c^{min}_a, c^{max}_a \in \mathbb{R}_+$. In other words, we prove that for any $a \in A$ and $\overline{c}_a \in [c^{min}_a, c^{max}_a]$ there exists a marginal probability distribution $\mathbb{Q}_a \in \widetilde{\mathcal{Q}}_a$ such that its expectation satisfies $\mathbb{E}_{\mathbb{Q}_a}\{c_a\} = \overline{c}_a$. The next result holds.
\begin{lemma} \label{lemma3}
Fix $a \in A$ and consider the set of marginal distributions $\widetilde{\mathcal{Q}}_a$ given by (\ref{set of marginal distributions}). Define
\begin{subequations}
\begin{align}
& c^{min}_a := \min_{\mathbb{Q}_a \in \widetilde{\mathcal{Q}}_a} \mathbb{E}_{\mathbb{Q}_a}\{c_a\} \label{best case expected cost}\\
& c^{max}_a := \max_{\mathbb{Q}_a \in \widetilde{\mathcal{Q}}_a} \mathbb{E}_{\mathbb{Q}_a}\{c_a\} \label{worst case expected cost}
\end{align}
\end{subequations}
Then for any $\overline{c}_a \in [c^{min}_a, c^{max}_a]$ there exists a marginal probability distribution $\mathbb{Q}_a \in \widetilde{\mathcal{Q}}_a$ such that $\mathbb{E}_{\mathbb{Q}_a}\{c_a\} = \overline{c}_a$.
\begin{proof} We need to verify whether the following set of probability distributions is nonempty:

\begin{equation} \label{set of marginal distributions with expectation constraint}
\begin{gathered}
\Big\{\mathbb{Q}_a \in \mathcal{Q}_0(\mathbb{R}): \mathbb{Q}_a\{c_a \in [l^{(i)}_a, u^{(i)}_a]\} \in [\underline{q}^{(i)}_a, \overline{q}^{(i)}_a] \quad \forall i \in \mathcal{D}_a \\ \mathbb{E}_{\mathbb{Q}_a}\{c_a\} =~\overline{c}_a\Big\}
\end{gathered}
\end{equation}
In order to establish this fact, we construct the corresponding feasibility problem:
\allowdisplaybreaks\begin{subequations} \label{feasibility problem with expectation constraint}
\begin{align}
& \max 0 \\
\mbox{\upshape s.t. }
& \mu \in \mathcal{M}_+(\mathbb{R}) \\
& \int_{l_a}^{u_a} \mathbbm{1}\Big\{l^{(i)}_a \leq c_a \leq u^{(i)}_a\Big\} \mbox{ d}\mu(c_a) \leq \overline{q}^{(i)}_a \quad \forall i \in \mathcal{D}_a \\
& \int_{l_a}^{u_a} \mathbbm{1}\Big\{l^{(i)}_a \leq c_a \leq u^{(i)}_a\Big\} \mbox{ d}\mu(c_a) \geq \underline{q}^{(i)}_a \quad \forall i \in \mathcal{D}_a, \\
& \int_{l_a}^{u_a} c_a \mbox{ d}\mu(c_a) = \overline{c}_a \label{cons: feasibility problem with expectation constraint 3}
\end{align}
\end{subequations}

In the sequel, we fix $a \in A$. Following the proof of Lemma \ref{lemma2} we obtain a dual reformulation of the feasibility problem (\ref{feasibility problem with expectation constraint}), i.e.,
\allowdisplaybreaks\begin{subequations} \label{individual dual problem with expectation constraint}
\begin{align}
& \min_{\mathbf{k}_a,\mathbf{h}_a, \gamma} \sum_{i \in \mathcal{D}_a}(\overline{q}^{(i)}_a k_{ai} - \underline{q}^{(i)}_a h_{ai}) + \overline{c}_a \gamma \\
\mbox{\upshape s.t. }
& k_{ai} \geq 0, \mbox{ } h_{ai} \geq 0 \quad \forall i \in \mathcal{D}_a \label{cons: individual dual problem with expectation constraint 1}\\
& \min_{j \in \mathcal{W}_a}\Big\{\sum_{i \in \mathcal{D}_a(j)}(k_{ai} - h_{ai}) + \min\{\gamma L^{(j)}_a, \gamma U^{(j)}_a\}\Big\} \geq 0 \label{cons: individual dual problem with expectation constraint 2}
\end{align}
\end{subequations}
Here, $\gamma \in \mathbb{R}$ is a dual variable corresponding to the primal constraint (\ref{cons: feasibility problem with expectation constraint 3}). In particular, notice that the feasible region of (\ref{individual dual problem with expectation constraint}) is nonempty since the zero solution, i.e., $\mathbf{k}_a, \mathbf{h}_a = \mathbf{0}$ and $\gamma = 0$, is feasible.

Hence, by the strong duality; see \cite{Wiesemann2014, Shapiro2001}, either primal optimization problem (\ref{feasibility problem with expectation constraint}) is feasible, or there exists a dual feasible solution $\widetilde{\mathbf{k}}_a, \widetilde{\mathbf{h}}_a \in \mathbb{R}^{D_a}$, $\widetilde{\gamma} \in \mathbb{R}$ such that:
\begin{equation} \label{dual unbounded}
\sum_{i \in \mathcal{D}_a}(\overline{q}^{(i)}_a \widetilde{k}_{ai} - \underline{q}^{(i)}_a \widetilde{h}_{ai}) + \overline{c}_a\widetilde{\gamma} < 0
\end{equation}
In fact, inequality (\ref{dual unbounded}) indicates that the dual problem (\ref{individual dual problem with expectation constraint}) is unbounded as long as the dual feasible solution $(\widetilde{\mathbf{k}}_a, \widetilde{\mathbf{h}}_a, \widetilde{\gamma})^\top$ can be multiplied by any positive constant.
Therefore, it is sufficient to show that for any dual feasible solution $\mathbf{k}_a, \mathbf{h}_a \in \mathbb{R}^{D_a}$, $\gamma \in \mathbb{R}$ inequality (\ref{dual unbounded}) does not hold, that is, we have:
\begin{equation} \label{dual bounded}
\sum_{i \in \mathcal{D}_a}(\overline{q}^{(i)}_a k_{ai} - \underline{q}^{(i)}_a h_{ai}) + \overline{c}_a \gamma \geq 0
\end{equation}

First, assume that $\gamma = -1$. Then the feasible region of optimization problem (\ref{individual dual problem with expectation constraint}) coincides with the feasible region of dual problem (\ref{individual dual problem}); recall the proof of Lemma~\ref{lemma2}. Hence, for any dual feasible solution $(\mathbf{k}_a, \mathbf{h}_a, -1)^\top$ we have:
$$\sum_{i \in \mathcal{D}_a}(\overline{q}^{(i)}_a k_{ai} - \underline{q}^{(i)}_a h_{ai}) - \overline{c}_a \geq c^{max}_a - \overline{c}_a \geq 0$$
Here, the first inequality is implied by the strong duality, i.e., the optimal objective function values in (\ref{worst case expected cost}) and (\ref{individual dual problem}) coincide, whereas the second inequality stems from the assumption $\overline{c}_a \in [c^{min}_a, c^{max}_a]$. We conclude that (\ref{dual bounded}) holds.

It is rather straightforward to verify that a dual reformulation of the minimization problem
$$\min_{\mathbb{Q}_a \in \widetilde{\mathcal{Q}}_a} \mathbb{E}_{\mathbb{Q}_a}\{c_a\}$$
is given by:
\allowdisplaybreaks\begin{subequations} \label{individual dual problem for min}
\begin{align}
& \max_{\mathbf{k}_a,\mathbf{h}_a} \sum_{i \in \mathcal{D}_a} -(\overline{q}^{(i)}_a k_{ai} - \underline{q}^{(i)}_a h_{ai})\\
\mbox{\upshape s.t. }
& k_{ai} \geq 0, \mbox{ } h_{ai} \geq 0 \quad \forall i \in \mathcal{D}_a \\
& \min_{j \in \mathcal{W}_a}\Big\{\sum_{i \in \mathcal{D}_a(j)}(k_{ai} - h_{ai}) + L^{(j)}_a\Big\} \geq 0
\end{align}
\end{subequations}
Thus, if $\gamma = 1$, then the feasible regions of (\ref{individual dual problem with expectation constraint}) and (\ref{individual dual problem for min}) coincide. Analogously, for any dual feasible solution $(\mathbf{k}_a, \mathbf{h}_a, 1)^\top$ we have:
$$\sum_{i \in \mathcal{D}_a}(\overline{q}^{(i)}_a k_{ai} - \underline{q}^{(i)}_a h_{ai}) + \overline{c}_a \geq -c^{min}_a + \overline{c}_a \geq 0$$
and thus, (\ref{dual bounded}) holds.
Furthermore, for any $\gamma \neq \pm 1$ the result is induced by scaling of the parameters, i.e., by introducing new endpoints given by:
\begin{eqnarray}
\widetilde{L}^{(j)}_a := \gamma L^{(j)}_a, \mbox{ if } \gamma \geq 0 \nonumber\\
\widetilde{U}^{(j)}_a := -\gamma U^{(j)}_a, \mbox{ if } \gamma < 0 \nonumber
\end{eqnarray}
This observation concludes the proof.
\end{proof}
\end{lemma}

Lemma \ref{lemma3} provides an intuition behind the robust reformulation of DRSPP (\ref{DRO formulation}). Indeed, quantile constraints (\ref{eq: quantile constraints}) can be replaced by the constructed linear expectation constraints (\ref{new linear expectation constraints}) with respect to each $a \in A$. In particular, the values of $c^{min}_a$ and $c^{max}_a$ are derived by solving the individual moment problems (\ref{best case expected cost}) and (\ref{worst case expected cost}), respectively. These problems can be solved in polynomial time via their dual reformulations; see the proofs of Theorem \ref{theorem1} and Lemma \ref{lemma3}. As a consequence, the next result provides a reformulation of (\ref{DRO formulation}) as an instance of the robust shortest path problem with polyhedral uncertainty.
\begin{theorem} \label{theorem2} Let
\begin{equation} \nonumber
\mathcal{S}: = \{\overline{\mathbf{c}} \in \mathbb{R}^{|A|}:\mbox{ } \mathbf{c}^{min} \leq \overline{\mathbf{c}} \leq \mathbf{c}^{max}; \quad \mathbf{B}\overline{\mathbf{c}} \leq \mathbf{b}\}
\end{equation}
Assume that $\mathbf{c}^{min} = \{c^{min}_a, a \in A\}$ and $\mathbf{c}^{max} = \{c^{max}_a, a \in A\}$ are given by (\ref{best case expected cost}) and (\ref{worst case expected cost}), respectively. Then the distributionally robust shortest path problem (\ref{DRO formulation}) is equivalent to the following robust shortest path problem with polyhedral uncertainty:
\begin{equation} \label{robust formulation} \nonumber \tag{\textbf{F2}}
\min_{\mathbf{y} \in Y} \max_{\overline{\mathbf{c}} \in \mathcal{S}} \sum_{a \in A} \overline{c}_a y_a
\end{equation}
\begin{proof} The result follows from Lemma \ref{lemma3} by setting $\mathbb{E}_{\mathbb{Q}_a}\{c_a\} = \overline{c}_a$, $a \in A$.
\end{proof}
\end{theorem}
Importantly, under Assumptions \textbf{A1} and \textbf{A2} Theorem \ref{theorem2} can be applied to any combinatorial optimization problem with distributional constraints (\ref{eq: linear expectation constraints}), (\ref{eq: support constraint}) and (\ref{eq: quantile constraints}). In particular, bounds (\ref{best case expected cost}) and (\ref{worst case expected cost}) for the cost vector $\mathbf{c}$ can be computed in polynomial time. 
Next, we explore complexity of RSPP (\ref{robust formulation}) and briefly discuss the associated solution techniques.
\subsection{Complexity and solution approach} \label{subsec: complexity}
We deduce that RSPP~(\ref{robust formulation}) is $NP$-hard by leveraging an equivalent robust formulation with a finite number of scenarios. The latter problem is known to be $NP$-hard even for a restricted class of networks, i.e., for layered graphs of width~$2$ \cite{Yu1998}. Regarding the solution techniques we employ the structure of the second-level optimization problem in (\ref{robust formulation}) to derive a single-level linear MIP reformulation.

The complexity results of this section are similar to the results discussed in \cite{Dokka2017}. However, we reiterate some basic ideas to preserve consistency and completeness of the manuscript.

First, observe that RSPP (\ref{robust formulation}) can be introduced as follows:
\allowdisplaybreaks\begin{subequations} \label{robust formulation complexity}
\begin{align}
& \min z\\
\mbox{\upshape s.t. }
& z \geq \overline{\mathbf{c}}^\top \mathbf{y} \quad \forall \overline{\mathbf{c}} \in \mathcal{S} \label{eq: robust constraint complexity}\\
& \mathbf{y} \in Y
\end{align}
\end{subequations}
Observe that the polyhedron $\mathcal{S}$ is bounded due to the added linear expectation constraints (\ref{new linear expectation constraints}). Let $\mathbf{c}^{(1)}, \ldots, \mathbf{c}^{(m)}$ be a vertex representation \cite{Conforti2014} of the polyhedral uncertainty set $\mathcal{S}$, i.e.,
$$\mathcal{S} = Conv(\mathbf{c}^{(1)}, \ldots, \mathbf{c}^{(m)})$$
Specifically, $Conv(\mathbf{c}^{(1)}, \ldots, \mathbf{c}^{(m)})$ denotes a convex hull of $\mathbf{c}^{(1)}, \ldots, \mathbf{c}^{(m)}$. Then RSPP (\ref{robust formulation}) is equivalent to the robust shortest path problem with $m$ discrete scenarios given by:
\allowdisplaybreaks\begin{subequations} \label{scenario formulation}
\begin{align}
& \min z\\
\mbox{\upshape s.t. }
& z \geq \mathbf{c}^{(i) \top} \mathbf{y} \quad \forall i \in \{1, \ldots, m\} \label{eq: robust constraint scenario}\\
& \mathbf{y} \in Y
\end{align}
\end{subequations}

More precisely, any constraint in (\ref{eq: robust constraint complexity}) can be represented as a convex combination of the constraints in (\ref{eq: robust constraint scenario}). Vice versa, constraints (\ref{eq: robust constraint scenario}) are contained in (\ref{eq: robust constraint complexity}).

Optimization problem in the form of (\ref{scenario formulation}) is known to be $NP$-hard even for two scenarios and for layered networks of width $2$ \cite{Yu1998}. We conclude that RSPP (\ref{robust formulation}) is $NP$-hard for $m = 2$ since we can reduce the two-scenario case by choosing $\mathcal{S}$ as a convex hull of the two scenarios. As a remark, the complexity class of (\ref{robust formulation}) for a general polyhedral uncertainty set is undefined as the number of facets of $\mathcal{S}$ might be exponential in $m$. We refer the interested reader to Theorem 4 in \cite{Buchheim2018} for a more comprehensive discussion. 


Next, we provide a linear mixed-integer programming reformulation of DRSPP (\ref{DRO formulation}) by dualizing the lower-level optimization problem in the robust formulation (\ref{robust formulation}); see, e.g., \cite{Audet1997}. The following result concludes our theoretical analysis.
\begin{theorem} \label{theorem3}
Distributionally robust shortest path problem (\ref{DRO formulation}) admits a mixed-integer programming reformulation:
\allowdisplaybreaks\begin{subequations}
\begin{align} \label{mip formulation} \nonumber \tag{\textbf{F3}}
& \min_{\mathbf{y}, \boldsymbol{\lambda}, \boldsymbol{\mu}, \boldsymbol{\nu}} \mathbf{b}^\top \boldsymbol{\lambda} + (\mathbf{c}^{max})^\top \boldsymbol{\nu} - (\mathbf{c}^{min})^\top \boldsymbol{\mu} \nonumber\\
\mbox{\upshape s.t. }
& \boldsymbol{\lambda}, \boldsymbol{\mu}, \boldsymbol{\nu} \geq 0 \nonumber\\
& -\mathbf{y} + \mathbf{B}^\top \boldsymbol{\lambda} + \boldsymbol{\nu} = \boldsymbol{\mu} \nonumber\\
& \mathbf{y} \in Y \nonumber,
\end{align}
\end{subequations}
where $\mathbf{c}^{min}$ and $\mathbf{c}^{max}$ are given by (\ref{best case expected cost}) and (\ref{worst case expected cost}), respectively.
\begin{proof}
Consider the robust reformulation (\ref{robust formulation}) of DRSPP (\ref{DRO formulation}). Notice that for fixed $\mathbf{y} \in Y$ the lower-level maximization problem
\begin{equation} \label{lower-level linear program}
\max_{\overline{\mathbf{c}} \in \mathcal{S}} \mathbf{c}^\top \mathbf{y}
\end{equation}
is a linear program. Hence, RSPP (\ref{robust formulation}) can be viewed as a single-level MIP problem by dualizing~(\ref{lower-level linear program}). This observation results in formulation (\ref{mip formulation}) and concludes the proof.
\end{proof}
\end{theorem}

MIP problem (\ref{mip formulation}), in turn, can be tackled using off-the-shelf mixed-integer programming software. Numerical results of the next section allude that MIP problem (\ref{mip formulation}) can be solved reasonably fast even for large-scale problem instances.

\section{Computational study} \label{sec: comp study}
We test DRSPP (\ref{DRO formulation}) on a class of synthetic randomly generated instances. We also apply the methodology of Section \ref{subsec: modeling the ambiguity set} for modeling the ambiguity set and analyse the model's out-of-sample performance. Specifically, our approach is compared with some other robust and distributionally robust optimization techniques in terms of the \textit{nominal expected loss}. Additionally, we analyze the solution times of MIP formulation~(\ref{mip formulation}).

The remainder of this section is organized as follows. Section \ref{subsec: benchmark} describes the alternative approaches. In Section \ref{subsec: test instances} we discuss how to generate the test instances. In Sections \ref{subsec: results} and \ref{subsec: summary} we provide a brief discussion of our numerical results and conclusions, respectively.

\subsection{Benchmark approaches} \label{subsec: benchmark}
 With respect to robust solution techniques, we consider the budget constrained approach of Bertsimas~et~al.~\cite{Bertsimas2003}. The authors introduce a parameter $\Gamma \in \{0, \ldots, |A|\}$, which corresponds to the maximal number of cost coefficients that deviate from their nominal values; the rest coefficients are set to their lower bounds. That is, we solve the following optimization problem:
\begin{equation} \label{budget constrained formulation} \nonumber \tag{\textbf{R$_0$}}
\min_{\mathbf{y} \in Y}\Big\{\mathbf{l}^\top \mathbf{y} + \max_{\{A' \subseteq A: |A'| \leq \Gamma\}}\sum_{a \in A'}(u_a - l_a)y_a\Big\}
\end{equation}

For fixed $\Gamma \in \{0, \ldots, |A|\}$ optimization problem (\ref{budget constrained formulation}) can be solved by considering $|A| + 1$ deterministic shortest path problems; we refer to \cite{Bertsimas2003} for more details. The choice of $\Gamma$ adjusts the robustness of the proposed method against the level of conservatism of the solution.

Also, we compare our model with a specialized distributionally robust model. We adopt a marginal moment ambiguity set from \cite{Bertsimas2004} by leveraging support constraints as well as an upper bound for the first- and second-order moments. That is, the ambiguity set has the form:
\begin{equation} \label{ambiguity set 2}
\begin{gathered}
\mathcal{Q}_2 := \Big\{\mathbb{Q} \in \mathcal{Q}_0(\mathbb{R}^{|A|}) \mbox{: } {\mathbb{Q}_a} \{l_a \leq c_a \leq u_a\} = 1,\\
\qquad \quad \mathbb{E}_{\mathbb{Q}_a} \{c_a\} \leq \widehat{\mu}_a, \mbox{ }
 \mathbb{E}_{\mathbb{Q}_a} \{c_a^2\} \leq \widehat{\sigma}_a^2, \mbox{ } a \in A
\Big\}
\end{gathered}
\end{equation}
and the distributionally robust problem is given by:
\begin{align} \label{DRO formulation 2} \nonumber \tag{\textbf{DR$_0$}}
\min_{\mathbf{y} \in Y} \max_{\mathbb{Q} \in \mathcal{Q}_2} \mathbb{E}_{\mathbb{Q}}\{\mathbf{c}^\top \mathbf{y}\}
\end{align}

For simplicity we assume that $\widehat{\mu}_a \leq u_a$ and $\widehat{\sigma}_a \leq u_a$ for each $a \in A$; otherwise the corresponding moment constraint is redundant. Similar to the problem without linear expectation constraints (\ref{DRO formulation without EC}), optimization problem (\ref{DRO formulation 2}) is partitioned into $|A|$ individual moment problems and a single deterministic shortest path problem.

It is rather straightforward to verify that the optimal expected costs in (\ref{DRO formulation 2}) are given by $\min\{\widehat{\mu}_a; \widehat{\sigma}_a\}$, $a \in A$. Furthermore, for any $a \in A$ the worst-case distribution is a discrete distribution, which takes value $\min\{\widehat{\mu}_a; \widehat{\sigma}_a\}$ with probability $1$. In particular, this fact follows from the nonnegativity of variance, i.e.,
\begin{equation} \nonumber
(\mathbb{E}_{\mathbb{Q}_a} \{c_a\})^2 \leq \mathbb{E}_{\mathbb{Q}_a} \{c_a^2\} \leq \widehat{\sigma}_a^2
\end{equation}
Next, we discuss generation of the test instances.

\subsection{Test instances } \label{subsec: test instances}
\begin{figure}
\centering
\begin{tikzpicture}[scale=0.75,transform shape]
\Vertex[x=0,y=0]{1}
\Vertex[x=3,y=3]{2}
\Vertex[x=3,y=0]{3}
\Vertex[x=3,y=-3]{4}
\Vertex[x=6,y=3]{5}
\Vertex[x=6,y=0]{6}
\Vertex[x=6,y=-3]{7}
\Vertex[x=9,y=3]{8}
\Vertex[x=9,y=0]{9}
\Vertex[x=9,y=-3]{10}
\Vertex[x=12,y=0]{11}
\tikzstyle{LabelStyle}=[fill=white,sloped]
\tikzstyle{EdgeStyle}=[post]
\Edge(1)(2)
\Edge(1)(3)
\Edge(1)(4)
\Edge(8)(11)
\Edge(9)(11)
\Edge(10)(11)
\Edge(2)(5)
\Edge(2)(6)
\Edge(2)(7)
\Edge(3)(5)
\Edge(3)(6)
\Edge(3)(7)
\Edge(4)(5)
\Edge(4)(6)
\Edge(4)(7)
\Edge(5)(8)
\Edge(5)(9)
\Edge(5)(10)
\Edge(6)(8)
\Edge(6)(9)
\Edge(6)(10)
\Edge(7)(8)
\Edge(7)(9)
\Edge(7)(10)
\end{tikzpicture}
\caption{\footnotesize A fully-connected layered graph with $v = 3$ intermediate layers and $r_i = 3$, $i \in \{1,\ldots, v\}$ nodes at each layer.}
\label{fig: layered graph}
\end{figure}
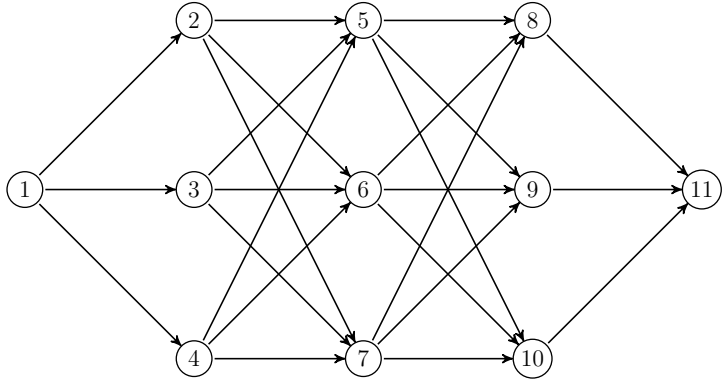

In our experiments we consider \textit{a fully-connected layered graph} with $v$ intermediate layers and $r_i$ nodes at each layer $i \in \{1, \ldots, v\}$. The first and the last layer consist of unique nodes, which are the source and the destination nodes, respectively, i.e., with some abuse of notation let $r_0 = r_{v + 1} = 1$. For example, a network with $v = 3$ and $r_i = 3$, $i \in \{1,2,3\}$, is depicted in Figure \ref{fig: layered graph}.

\subsubsection{Construction of the nominal distribution}
We construct a \textit{nominal distribution} $\mathbb{Q}^0$ of the cost vector $\mathbf{c}$ as a product of the corresponding marginal distributions $\mathbb{Q}^0_a \in \mathcal{Q}_0(\mathbb{R})$, $a \in A$. For each $a \in A$ we construct the support by setting $l_a$ uniformly distributed on $[0, 100]$, i.e., $l_a \sim \mathcal{U}(0, 100)$, and $u_a := l_a + \Delta_a$, where $\Delta_a \sim \mathcal{U}(0, 100)$.

Assume that for fixed $a \in A$ the arc cost $c_a$ is governed by a generalized \textit{beta distribution} with parameters $\alpha_a, \beta_a \in \mathbb{R}_+$ and the support given by $[l_a, u_a]$. Denote by $m_a$ and $\sigma_a$ its mean and variance, respectively. By standard calculations \cite{Gupta2004} observe that:
\begin{equation} \label{eq: parameters of beta distribution}
\begin{gathered}
\alpha_a = \frac{\widetilde{m}_a^2(1 - \widetilde{m}_a)}{\widetilde{\sigma}_a} - \widetilde{m}_a\\
\beta_a = \alpha_a(\frac{1}{\widetilde{m}_a} - 1),
\end{gathered}
\end{equation}
where $\widetilde{m}_a := (m_a - l_a)/(u_a - l_a)$ and $\widetilde{\sigma}_a := \sigma_a/(u_a - l_a)^2$ are the normalized mean and variance.

Now, to construct a beta distribution for each particular arc $a \in A$ we set $\widetilde{\sigma}_a = 1/64$ and
\begin{equation} \label{eq: nominal expected cost} \nonumber
\widetilde{m}_a \sim \mathcal{U}(\frac{1}{2}(1 - \sqrt{1 - 4\widetilde{\sigma}_a}), \frac{1}{2}(1 + \sqrt{1 - 4\widetilde{\sigma_a}})), \quad \widetilde{m}_a \neq \frac{1}{2}(1 \pm \sqrt{1 - 4\widetilde{\sigma}_a})
\end{equation}
Indeed, the abovementioned conditions stipulate that a beta distribution with the parameters $\alpha_a$, $\beta_a$ defined by (\ref{eq: parameters of beta distribution}) exists, i.e., $\alpha_a, \beta_a > 0$.

\subsubsection{Construction of the distributional constraints}
\begin{table}

\footnotesize
  \centering
  \onehalfspacing
    \begin{tabular}{c c c}
    \hline
    Parameter & Definition & Value(s) \\
    \hline
    $n_0 \in \mathbb{Z}_+$ & number of samples for each $a \in A$ & 100 \\
    $n_1 \in \mathbb{Z}_+$ & number of subintervals for each $a \in A$ & $n_1 \in \{1, 2, 3, 4\}$ \\
    $\eta_0 \in (0,1)$ & probability of violation for the ambiguity set $\mathcal{Q}$ & 0.05 \\
    $\eta \in (0,1)$ & probability of violation for each individual constraint & $\frac{\eta_0}{n_1|A| + D_0}$  \\
    $\kappa \in (0,1)$ & relative width, $\frac{u^{(i)}_a - l^{(i)}_a}{u_a - l_a}$, for each $i \in \mathcal{D}_a \setminus \{1\}$ and $a \in A$ & $\kappa \in \{0.2, 0.4, 0.6, 0.8\}$ \\
    \hline
    \end{tabular}
  \caption{Parameters of the ambiguity set. Recall that $\mathcal{D}_a = \{1, \ldots, D_a\}$ are the indices of baseline subintervals for each $a \in A$; $D_0$ is the number of linear expectation constraints.}
  \label{tab: parameters}
\end{table}

We apply the methodology of Section \ref{subsec: modeling the ambiguity set} to construct our distributional constraints from historical data. The summary of our parameter's settings is provided in Table \ref{tab: parameters}. Initially, for each arc $a \in A$ we generate $n_0$ random i.i.d. observations of $c_a$, which are referred to as $\xi^{(a)}_i$, $i \in \{1, \ldots, n_0\}$.

At first, we use (\ref{eq: to model quantiles}) to generate quantile constraints in the form of (\ref{eq: quantile constraints}). Given relative width $\kappa \in (0,1)$ we generate the subintervals $[l^{(i)}_a, u^{(i)}_a], i \in \mathcal{D}_a \setminus \{1\}$, by setting
$$l^{(i)}_a \sim \mathcal{U}(l_a, u_a - \kappa(u_a - l_a))$$
and $u^{(i)}_a := l^{(i)}_a + \kappa(u_a - l_a)$; see Table \ref{tab: parameters}. Then we introduce a probability of violation, $\eta \in (0,1)$, for each individual constraint and follow the procedure described in Section \ref{subsec: modeling the ambiguity set}.



Next, we model linear expectation constraints (\ref{eq: linear expectation constraints}) in the following way. Initially, we construct a subset of ``near-optimal'' paths $\widetilde{\mathcal{P}} \subseteq \mathcal{P}_{st}(G)$ with respect to the worst-case expected costs $c^{max}_a$, $a \in A$; recall Lemma \ref{lemma1}. Specifically, let $\widetilde{P}^*$ be an optimal path induced by (\ref{DRO formulation without EC}). For each arc $a \in \widetilde{P}^*$ we remove this arc from the network $G$ and seek the shortest path $\widetilde{P}^*_a$ in the resulting network $G[A \setminus a]$.
As a result, assume that the set $\widetilde{\mathcal{P}}$ is comprised of the path $\widetilde{P}^*$ and the newly constructed paths $\widetilde{P}^*_a$ for each $a \in \widetilde{P}^*$.
For each path $P \in \widetilde{\mathcal{P}}$ we introduce the associated linear expectation constraints (\ref{eq: to model linear expectation constraints}) subject to probability of violation $\eta$.

\begin{figure}
\centering
\begin{tikzpicture}[scale=0.75,transform shape]
\Vertex[x=0,y=0]{1}
\Vertex[x=6,y=2.5]{2}
\Vertex[x=6,y=-2.5]{3}
\Vertex[x=12,y=0]{4}
\tikzstyle{LabelStyle}=[fill=white,sloped]
\tikzstyle{EdgeStyle}=[post]
\tikzstyle{EdgeStyle}=[post, bend left, thick, double = black,
                                 double distance = 0.5pt]
\Edge[label=$37$](1)(2)
\Edge[label=$101$](2)(4)
\tikzstyle{EdgeStyle}=[post, bend left]
\Edge[label=$100$](2)(3)
\tikzstyle{EdgeStyle}=[post, bend right]
\Edge[label=$100$](1)(3)
\Edge[label=$100$](3)(4)
\end{tikzpicture}
\caption{\footnotesize The network used in Example \ref{ex: Example 1}. The worst-case expected costs with regard to formulation (\ref{DRO formulation without EC}), i.e., $c^{max}_a$, are depicted inside each arc $a \in A$. The shortest path $\widetilde{P}^*$ is in bold.}
\label{fig: Example 1 continued}
\end{figure}
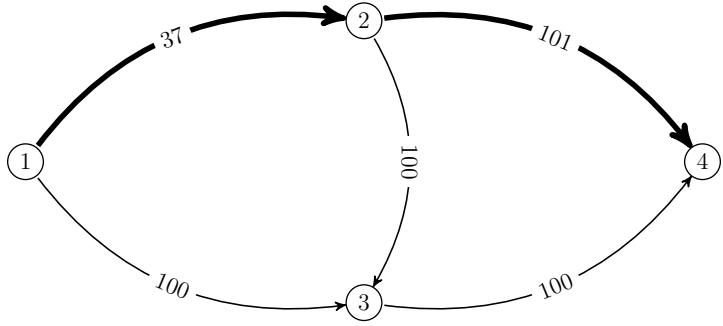

\begin{repexample}{ex: Example 1}[Continued] \upshape Consider the network used in Example \ref{ex: Example 1}; see Figure \ref{fig: Example 1 continued}. The shortest path $\widetilde{P}^*$ with regard to $\mathbf{c}^{max}$ is given by $\{1 \rightarrow 2 \rightarrow 4\}$. If we remove either arc $(1, 2)$, or arc $(2, 4)$ from the network, then the shortest path in the resulting network is given by $\{1 \rightarrow 3 \rightarrow 4\}$. Hence, the set $\widetilde{\mathcal{P}}$ consists of two distinct paths $\{1 \rightarrow 2 \rightarrow 4\}$ and $\{1 \rightarrow 3 \rightarrow 4\}$. \vspace{-9mm}\flushright$\square$
\end{repexample}

Observe that by construction Assumption \textbf{A1} holds since the discrete uniform distribution at the generated samples satisfies both quantile and linear expectation constraints. Furthermore, Assumption \textbf{A2} is satisfied due to randomness in the choice of baseline subintervals. As outlined in Section \ref{subsec: modeling the ambiguity set} we use Bonferroni's inequality to provide a required confidence level for the ambiguity set $\mathcal{Q}$. More precisely, we denote by $\eta_0 \in (0,1)$ the probability that at least one distributional constraint contained in $\mathcal{Q}$ is violated; see Table \ref{tab: parameters}.

Finally, we briefly discuss construction of the ambiguity set (\ref{ambiguity set 2}) for the moment-based formulation (\ref{DRO formulation 2}). Recall that to construct (\ref{eq: to model quantiles}) we only need to verify whether random samples $\xi^{(a)}_j$, $j \in \{1, \ldots, n_0\}$, belong to a prescribed subinterval $[l^{(i)}_a, u^{(i)}_a]$, $i \in \mathcal{D}_a$, $a \in A$, or not. Hence, to provide a fair comparison we replace pointwise observations $\xi^{(a)}_j$ (which can be unavailable to the decision-maker) with their maximal possible values. Then we estimate $\widehat{\mu}_a$ and $\widehat{\sigma}_a^2$ for each $a \in A$; we refer to \ref{sec: app} for further details.


\subsubsection{Computational settings}
All experiments are performed on a PC with \textit{CPU i5-7200U} and \textit{RAM 6 GB}. MIP problem (\ref{mip formulation}) is solved in Java with \textit{CPLEX 12.7.1}. The deterministic shortest path problems are solved with Dijkstra's algorithm \cite{Dijkstra1959}. Additionally, we set $v = 20$, $r_i = 10$, $i \in \{1,\ldots, 20\}$, $n_0 = 100$ and $\eta_0 = 0.05$ for the distributionally robust formulations (\ref{DRO formulation}) and (\ref{DRO formulation 2}). We consider several possible values of the relative width, $\kappa$, and the number of subintervals, $n_1$, that is, $\kappa \in \{0.2, 0.4, 0.6, 0.8\}$ and $n_1 \in \{1, 2, 3, 4\}$. Eventually, let $\Gamma \in \{0, 7, 14, 21\}$ in the budget constrained formulation~(\ref{budget constrained formulation}).

\subsection{Results and discussion} \label{subsec: results}
\subsubsection{Comparison with formulations (\ref{budget constrained formulation}) and (\ref{DRO formulation 2})}
We compare formulations (\ref{DRO formulation}), (\ref{DRO formulation without EC}) with the robust formulation (\ref{budget constrained formulation}) and distributionally robust formulation (\ref{DRO formulation 2}) in terms of their out-of-sample performance. In other words, for any path $\widehat{P} \in \mathcal{P}_{st}(G)$ we introduce a \textit{relative expected loss}
\begin{equation} \label{relative expected loss}
\rho^0(\widehat{P}) := \frac{\mathbb{E}_{\mathbb{Q}^0} {\sum_{a \in \widehat{P}} c_a}}{\min_{P \in \mathcal{P}_{st}(G)} \{\mathbb{E}_{\mathbb{Q}^0} {\sum_{a \in P} c_a}\}} = \frac{\sum_{a \in \widehat{P}} m_a}{\min_{P \in \mathcal{P}_{st}(G)} \{\sum_{a \in P} m_a\}},
\end{equation}
which enables to evaluate the quality of robust (distributionally robust) solutions. Specifically, (\ref{relative expected loss}) reflects the ratio of the nominal expected loss incurred by the decision-maker to the optimal expected loss in the full-information setting, i.e., when the expected values $m_a$, $a \in A$, are known in prior. In particular, observe that $\rho^0(\widehat{P}) \geq 1$.

Let $\kappa = 0.6$ and $n_1 = 4$. In Table \ref{tab: comparison} we report the average relative loss (\ref{relative expected loss}) incurred by the decision-maker across 100 random network instances.

\begin{table}
\footnotesize

  \centering
  \onehalfspacing
    \begin{tabular}{c c c}
    \hline
    Solution approach & Relative expected loss \\
    \hline
    DRSPP (\ref{DRO formulation without EC}) & 1.29 (0.12)\\
    DRSPP (\ref{DRO formulation}) & 1.22 (0.10) \\
    DRSPP (\ref{DRO formulation 2}) & 1.59 (0.15)\\
    RSPP (\ref{budget constrained formulation}) with $\Gamma = 0$ & 2.23 (0.27) \\
    RSPP (\ref{budget constrained formulation}) with $\Gamma = 7$ & 1.86 (0.23) \\
    RSPP (\ref{budget constrained formulation}) with $\Gamma = 14$ & 1.88 (0.22)  \\
    RSPP (\ref{budget constrained formulation}) with $\Gamma = 21$ & 2.19 (0.20) \\
    \hline
    \end{tabular}
  \caption{Let $\kappa = 0.6$ and $n_1 = 4$. We report the average relative losses (\ref{relative expected loss}) and standard deviations (in brackets) across 100 random instances.}
  \label{tab: comparison}
\end{table}


Going back to the discussion in Section \ref{sec: intro} observe that robust formulation (\ref{budget constrained formulation}) does not exploit any distributional information and thus, provides overly conservative solutions; see lines 3-7 of Table~\ref{tab: comparison}. Then in lines 1 and 3 we ensure that our formulation of DRSPP (\ref{DRO formulation without EC}) outperforms the moment-based formulation (\ref{DRO formulation 2}) in terms of the relative expected loss. In particular, both formulations (\ref{DRO formulation without EC}) and (\ref{DRO formulation 2}) assume no correlation information between the travel times. However, it is rather intuitive that the out-of-sample performance can be improved by utilizing linear expectation constraints; see line 2 in Table \ref{tab: comparison}. 



Next, we discuss dependence of the out-of-sample performance on the relative width, $\kappa$, and number of subintervals, $n_1$.


\subsubsection{Dependence on $\kappa$ and $n_1$}
Consider the formulation without linear expectation constraints (\ref{DRO formulation without EC}) and moment-based formulation (\ref{DRO formulation 2}) with $\eta_0 = 0.05$. We explore how parameters $\kappa$ and $n_1$ affect the quality of distributionally robust solutions in terms of the relative expected loss (\ref{relative expected loss}). That is, for $\kappa = 0.6$ we assume that $n_1 \in \{1, 2, 3, 4\}$; for $n_1 = 4$ we assume that $\kappa \in \{0.2, 0.4, 0.6, 0.8\}$. The results are reported in Tables \ref{tab: n1} and \ref{tab: kappa}, respectively. 

Naturally, with the increase of $n_1$ the relative expected loss decreases; see Table \ref{tab: n1}. In fact, despite the probability of violation for individual constraints increases; recall Table \ref{tab: parameters}, introducing new subintervals refines both formulations (\ref{DRO formulation without EC}) and (\ref{DRO formulation 2}). In particular, we append an additional quantile constraint to the ambiguity set (\ref{ambiguity set}) and specify the estimates of the first- and second-order moments in the ambiguity set (\ref{ambiguity set 2}).

\begin{table}
\footnotesize

  \centering
  \onehalfspacing
    \begin{tabular}{c c c c c c c c }
    \hline
    Solution approach & $n_1 = 1$ & $n_1 = 2$ & $n_1 = 3$ & $n_1 = 4$ \\
    \hline
    DRSPP (\ref{DRO formulation without EC}) & 1.59 (0.15) & 1.40 (0.12) & 1.31 (0.11) & 1.29 (0.12)\\
    DRSPP (\ref{DRO formulation 2}) & 1.82 (0.18) & 1.65 (0.16) & 1.61 (0.14) & 1.59 (0.15) \\
    \hline
    \end{tabular}
  \caption{ The relative expected losses (in average) and standard deviations (in brackets) across 100 random instances for fixed $\kappa = 0.6$ and $n_1 \in \{1, 2, 3, 4\}$. }
\label{tab: n1}
\end{table}

\begin{table}
\footnotesize

  \centering
  \onehalfspacing
    \begin{tabular}{c c c c c c c c }
    \hline
    Solution approach & $\kappa = 0.2$ & $\kappa = 0.4$ & $\kappa = 0.6$ & $\kappa = 0.8$ \\
    \hline
    DRSPP (\ref{DRO formulation without EC}) & 1.87 (0.17) & 1.39 (0.12) & 1.29 (0.12) & 1.32 (0.11)\\
    DRSPP (\ref{DRO formulation 2}) & 2.02 (0.19) & 1.75 (0.16) & 1.59 (0.15) & 1.62 (0.16)\\
     \hline
    \end{tabular}
  \caption{ The relative expected losses (in average) and standard deviations (in brackets) across 100 random instances for fixed $n_1 = 4$ and $\kappa \in \{0.2, 0.4, 0.6, 0.8\}$. }
\label{tab: kappa}
\end{table}

Further, with the increase of parameter $\kappa$ the relative expected loss decreases for $\kappa \leq 0.6$ and increases otherwise; see Table \ref{tab: kappa}. This observation for formulation (\ref{DRO formulation without EC}) can be justified as follows. For large values of $\kappa$ the subintervals in quantile constraints~(\ref{eq: to model quantiles}) contain a major part of the support $[l_a, u_a]$, $a \in A$. Hence, the distributionally robust approach provides ``near-robust'' decisions. In other words, for each $a \in A$ our estimate $c^{max}_a$ of the arc cost $c_a$ is sufficiently close to the upper bound $u_a$; this fact follows from the construction of the worst-case distribution with regard to (\ref{DRO formulation without EC}), see Lemma \ref{lemma2}. Alternatively, if $\kappa$ is small, then the probability that $c_a$ belongs to these subintervals is typically small, which results in ``near-robust'' decisions too. From Table \ref{tab: kappa} we conclude that intermediate values of $\kappa$, e.g., $\kappa = 0.6$, provide better solutions in terms of the relative expected loss~(\ref{relative expected loss}).

In a similar way, one may provide the intuition with regard to formulation (\ref{DRO formulation 2}). We also note that our formulation (\ref{DRO formulation without EC}) outperforms the moment-based formulation (\ref{DRO formulation 2}) across the considered values of parameters $\kappa$ and $n_1$.

\subsubsection{Running time}
In conclusion, we show that DRSPP (\ref{DRO formulation}) with both linear expectation (\ref{eq: to model linear expectation constraints}) and quantile constraints (\ref{eq: to model quantiles}) can be solved effectively using MIP solvers. Specifically, the solution procedure is divided into two stages:
\begin{enumerate}
\item[(\textit{i})] bounds $\mathbf{c}^{max}$ and $\mathbf{c}^{min}$ for the cost vector $\mathbf{c}$ are obtained by solving the individual moment problems (\ref{best case expected cost}) and (\ref{worst case expected cost}), respectively, for each $a \in A$;
\item[(\textit{ii})] the optimal solution of (\ref{DRO formulation}) is attained by solving the linear MIP problem (\ref{mip formulation}).
\end{enumerate}

Let $\eta_0 = 0.05$, $n_1 = 4$ and $\kappa = 0.6$. Assume that the number $v$ of intermediate layers is such that $v \in \{20, 40, \ldots, 100\}$ and there are $r_i = 10$, $i \in \{1, \ldots, v\}$, nodes at each layer. For each stage (\textit{i}) and (\textit{ii}) we report average and maximal running times over 100 randomly generated instances as a function of $v$; see Table~\ref{tab: running time}.

Recall that the cost bounds $c^{max}_a$ and $c^{min}_a$ for each $a \in A$ can be computed in polynomial time, while MIP formulation (\ref{mip formulation}) is $NP$-hard in general. Nevertheless, observe that the average solution times for the stage (\textit{ii}) are sufficiently small compared with the average solution times for the stage (\textit{i}). These observations can be explained by a specific construction of the linear expectation constraints (\ref{eq: to model linear expectation constraints}). In particular, the intuition is two-fold.
\begin{table}
\footnotesize

  \centering
  \onehalfspacing
    \begin{tabular}{c c c c c c c c }
    \hline
    & Running time & $v = 20$ & $v = 40$ & $v = 60$ & $v = 80$ & $v = 100$ \\
    \hline
    \multirow{2}{*}{Stage (\textit{i}) $\Bigg \{$} & Average running time (s) & 0.6 & 1.21 & 1.9 & 2.64 & 3.56\\
& Maximal running time (s) & 3.1 & 2.55 & 2.47 & 5.78 & 4.8\\
    \hline
   \multirow{2}{*}{Stage (\textit{ii}) $\Bigg \{$} & Average running time (s) & 0.06 & 0.12 & 0.17 & 0.3 & 0.5 \\
& Maximal running time (s) & 0.89 & 0.41 & 0.25 & 0.96 & 1.13\\
    \hline
    \end{tabular}
  \caption{Let $\eta_0 = 0.05$, $\kappa = 0.6$ and $n_1 = 4$. We report average and maximal running times for solving DRSPP (\ref{DRO formulation}) as a function of the network size, i.e., the number of intermediate layers $v$, over 100 random instances. }
  \label{tab: running time}
\end{table}

On the one hand, the number of linear expectation constraints is sufficiently small, that is, $D_0 = |\widetilde{\mathcal{P}}| = v + 2 \ll |A|$. Actually, guided by the discussion of Section \ref{subsec: modeling the ambiguity set} the construction of (\ref{eq: to model linear expectation constraints}) alludes that the decision-maker has a sufficient number of random observations of a total path cost. On the other hand, (\ref{eq: to model linear expectation constraints}) can provide initial feasible solutions for RSPP~(\ref{robust formulation}) with a reasonable linear programming relaxation quality. Namely, linear expectation constraints of the form (\ref{eq: to model linear expectation constraints}) bound the worst-case expected cost for paths $P \in \widetilde{\mathcal{P}}$ as well.

\subsection{Summary} \label{subsec: summary}
Summarizing the discussion above distributional constraints (\ref{eq: linear expectation constraints}) and (\ref{eq: quantile constraints}) provide a viable and flexible methodology to account for distributional uncertainty. The reported numerical results demonstrate the advantages of our approach against standard robust and distributionally robust optimization techniques. We conclude that the quality of distributionally robust solutions depends both on quality of information collected and a structure of the distributional constraints. In particular, with an appropriate choice of the parameters, distributionally robust solutions may provide a reasonable approximation of the optimal solution in the full-information setting. Moreover, we demonstrate that DRSPP (\ref{DRO formulation}) can be solved sufficiently fast using state-of-the-art mixed-integer programming solvers.

\section{Conclusion} \label{sec: conclusion}
In this paper we consider the shortest path problem, where the arc costs are governed by some probability distribution, which is itself subject to uncertainty. A distributionally robust version of the shortest path problem (DRSPP) is formulated, where the decision-maker attempts to minimize her worst-case expected loss over a family of candidate distributions that are compatible with the decision-maker's prior information.

The distributional family is formed by linear expectation constraints with respect to subsets of arcs and individual quantile constraints with respect to particular arcs. Our distributional constraints can be constructed from incomplete or partially observable data. For instance, one can employ observations of the cumulative cost for prescribed subsets of arcs and interval-censored observations for particular arcs. Then a family of candidate distributions is formed by leveraging standard measure concentration inequalities.

We propose equivalent robust and mixed-integer programming reformulations of DRSPP. In particular, the problem without linear expectation constraints is proved to be polynomially solvable. We demonstrate numerically that our approach is competitive against basic robust and   distributionally robust optimization techniques. Flexibility of the distributional constraints enables the decision-maker to collect distributional information and improve her solutions through multiple decision epochs. Furthermore, the proposed mixed-integer programming formulations can be solved effectively using state-of-the-art solvers.

Naturally, the theoretical results of this paper can be generalized to a class of polynomially solvable combinatorial optimization problems.
Furthermore, our ambiguity set can also be applied to conditional value at risk \cite{Rockafellar2000} and some other optimized certainty equivalent objective criteria \cite{Wiesemann2014, Hanasusanto2016}. To this end, one would need to apply the min-max theorem, dualization of the inner optimization problem and our construction of the elementary subintervals; recall the proof of Lemma 2. This approach also results in a mixed-integer programming reformulation of DRSPP, but after applying some auxiliary linearization techniques.  

As a future research direction, one can consider adapting the results of this paper for two- and multi-stage decision-making problems. For example, one can assume that the first part of a decision is picked before the realization of uncertainty, while the second part is determined after the uncertain costs are realized; see \cite{Goerigk2020}. An interesting concept of $K$-adaptability is discussed in \cite{Hanasusanto2016}, where a restricted number of preselected decisions is implemented at the second stage. Also, it may be possible to exploit the enhanced linear decision rules proposed in \cite{Bertsimas2018} to provide a fairly good approximation for some classes of two- or multi-stage distributionally robust formulations.


\textbf{Acknowledgments.} The article was prepared within the framework of the Basic Research Program at the National Research University Higher School of Economics (Sections 1-4) and funded by RFBR grant №20-37-90060 (Sections 5-6). The authors are thankful to the associate editor and three anonymous referees for their constructive comments that allowed us to greatly improve the paper.

\textbf{Conflict of interest:} The authors declare that they have no conflict of interest.

\bibliographystyle{ieeetr}
\bibliography{bibliography}

\begin{thebibliography}{10}

\bibitem{Dijkstra1959}
E.~W. Dijkstra, ``A note on two problems in connexion with graphs,'' {\em
  Numerische Mathematik}, vol.~1, no.~1, pp.~269--271, 1959.

\bibitem{Ahuja1988}
R.~K. Ahuja, T.~L. Magnanti, and J.~B. Orlin, {\em Network flows}.
\newblock Cambridge, Mass.: Alfred P. Sloan School of Management,
  Massachusetts, 1988.

\bibitem{Ford1958}
L.~R. Ford~Jr and D.~R. Fulkerson, ``A suggested computation for maximal
  multi-commodity network flows,'' {\em Management Science}, vol.~5, no.~1,
  pp.~97--101, 1958.

\bibitem{Johnson1978}
D.~S. Johnson, J.~K. Lenstra, and A.~R. Kan, ``The complexity of the network
  design problem,'' {\em Networks}, vol.~8, no.~4, pp.~279--285, 1978.

\bibitem{Israeli2002}
E.~Israeli and R.~K. Wood, ``Shortest-path network interdiction,'' {\em
  Networks}, vol.~40, no.~2, pp.~97--111, 2002.

\bibitem{Irnich2005}
S.~Irnich and G.~Desaulniers, ``Shortest path problems with resource
  constraints,'' in {\em Column generation} (G.~Desaulniers, J.~Desrosiers, and
  M.~M. Solomon, eds.), pp.~33--65, Springer, 2005.

\bibitem{Bellman1958}
R.~Bellman, ``On a routing problem,'' {\em Quarterly of Applied Mathematics},
  vol.~16, no.~1, pp.~87--90, 1958.

\bibitem{Bazaraa2011}
M.~S. Bazaraa, J.~J. Jarvis, and H.~D. Sherali, {\em Linear programming and
  network flows}.
\newblock John Wiley \& Sons, 2011.

\bibitem{Bertsimas2003}
D.~Bertsimas and M.~Sim, ``Robust discrete optimization and network flows,''
  {\em Mathematical Programming}, vol.~98, no.~1-3, pp.~49--71, 2003.

\bibitem{Mudchanatongsuk2008}
S.~Mudchanatongsuk, F.~Ord{\'o}{\~n}ez, and J.~Liu, ``Robust solutions for
  network design under transportation cost and demand uncertainty,'' {\em
  Journal of the Operational Research Society}, vol.~59, no.~5, pp.~652--662,
  2008.

\bibitem{Bental1998}
A.~Ben-Tal and A.~Nemirovski, ``Robust convex optimization,'' {\em Mathematics
  of Operations Research}, vol.~23, no.~4, pp.~769--805, 1998.

\bibitem{Bental2002}
A.~Ben-Tal and A.~Nemirovski, ``Robust optimization--methodology and
  applications,'' {\em Mathematical Programming}, vol.~92, no.~3, pp.~453--480,
  2002.

\bibitem{Bental2009}
A.~Ben-Tal, L.~El~Ghaoui, and A.~Nemirovski, {\em Robust optimization},
  vol.~28.
\newblock Princeton University Press, 2009.

\bibitem{Bertsimas2011}
D.~Bertsimas, D.~B. Brown, and C.~Caramanis, ``Theory and applications of
  robust optimization,'' {\em SIAM Review}, vol.~53, no.~3, pp.~464--501, 2011.

\bibitem{Ruszczynski2003}
A.~Ruszczy{\'n}ski and A.~Shapiro, ``Stochastic programming models,'' {\em
  Handbooks in Operations Research and Management Science}, vol.~10, pp.~1--64,
  2003.

\bibitem{Shapiro2009}
A.~Shapiro, D.~Dentcheva, and A.~Ruszczy{\'n}ski, {\em Lectures on stochastic
  programming: modeling and theory}.
\newblock SIAM, 2009.

\bibitem{Wiesemann2014}
W.~Wiesemann, D.~Kuhn, and M.~Sim, ``Distributionally robust convex
  optimization,'' {\em Operations Research}, vol.~62, no.~6, pp.~1358--1376,
  2014.

\bibitem{Zhang2017}
Y.~Zhang, S.~Song, Z.-J.~M. Shen, and C.~Wu, ``Robust shortest path problem
  with distributional uncertainty,'' {\em IEEE Transactions on Intelligent
  Transportation Systems}, vol.~19, no.~4, pp.~1080--1090, 2017.

\bibitem{Delage2010}
E.~Delage and Y.~Ye, ``Distributionally robust optimization under moment
  uncertainty with application to data-driven problems,'' {\em Operations
  research}, vol.~58, no.~3, pp.~595--612, 2010.

\bibitem{Goh2010}
J.~Goh and M.~Sim, ``Distributionally robust optimization and its tractable
  approximations,'' {\em Operations Research}, vol.~58, no.~4-part-1,
  pp.~902--917, 2010.

\bibitem{Bertsimas2018}
D.~Bertsimas, M.~Sim, and M.~Zhang, ``Adaptive distributionally robust
  optimization,'' {\em Management Science}, vol.~65, no.~2, pp.~604--618, 2018.

\bibitem{Kasperski2016}
A.~Kasperski and P.~Zieli{\'n}ski, ``Robust discrete optimization under
  discrete and interval uncertainty: A survey,'' in {\em Robustness analysis in
  decision aiding, optimization, and analytics} (M.~Doumpos1, C.~Zopounidis,
  and E.~Grigoroudis, eds.), pp.~113--143, Springer, 2016.

\bibitem{Dokka2017}
T.~Dokka and M.~Goerigk, ``An experimental comparison of uncertainty sets for
  robust shortest path problems,'' {\em arXiv preprint arXiv:1704.08470}, 2017.

\bibitem{Montemanni2004}
R.~Montemanni and L.~M. Gambardella, ``An exact algorithm for the robust
  shortest path problem with interval data,'' {\em Computers \& Operations
  Research}, vol.~31, no.~10, pp.~1667--1680, 2004.

\bibitem{Yu1998}
G.~Yu and J.~Yang, ``On the robust shortest path problem,'' {\em Computers \&
  Operations Research}, vol.~25, no.~6, pp.~457--468, 1998.

\bibitem{Zielinski2004}
P.~Zieli{\'n}ski, ``The computational complexity of the relative robust
  shortest path problem with interval data,'' {\em European Journal of
  Operational Research}, vol.~158, no.~3, pp.~570--576, 2004.

\bibitem{Hall1983}
R.~W. Hall, ``Travel outcome and performance: the effect of uncertainty on
  accessibility,'' {\em Transportation Research Part B: Methodological},
  vol.~17, no.~4, pp.~275--290, 1983.

\bibitem{Bertsekas1991}
D.~P. Bertsekas and J.~N. Tsitsiklis, ``An analysis of stochastic shortest path
  problems,'' {\em Mathematics of Operations Research}, vol.~16, no.~3,
  pp.~580--595, 1991.

\bibitem{Nie2009}
Y.~M. Nie and X.~Wu, ``Shortest path problem considering on-time arrival
  probability,'' {\em Transportation Research Part B: Methodological}, vol.~43,
  no.~6, pp.~597--613, 2009.

\bibitem{Chen2010}
A.~Chen and Z.~Zhou, ``The $\alpha$-reliable mean-excess traffic equilibrium
  model with stochastic travel times,'' {\em Transportation Research Part B:
  Methodological}, vol.~44, no.~4, pp.~493--513, 2010.

\bibitem{Pan2013}
Y.~Pan, L.~Sun, and M.~Ge, ``Finding reliable shortest path in stochastic
  time-dependent network,'' {\em Procedia-Social and Behavioral Sciences},
  vol.~96, pp.~451--460, 2013.

\bibitem{Shapiro2008}
A.~Shapiro, ``Stochastic programming approach to optimization under
  uncertainty,'' {\em Mathematical Programming}, vol.~112, no.~1, pp.~183--220,
  2008.

\bibitem{Verweij2003}
B.~Verweij, S.~Ahmed, A.~J. Kleywegt, G.~Nemhauser, and A.~Shapiro, ``The
  sample average approximation method applied to stochastic routing problems: a
  computational study,'' {\em Computational Optimization and Applications},
  vol.~24, no.~2-3, pp.~289--333, 2003.

\bibitem{Ng2011}
M.~Ng, W.~Szeto, and S.~T. Waller, ``Distribution-free travel time reliability
  assessment with probability inequalities,'' {\em Transportation Research Part
  B: Methodological}, vol.~45, no.~6, pp.~852--866, 2011.

\bibitem{Rinehart2011}
M.~Rinehart and M.~A. Dahleh, ``The value of side information in shortest path
  optimization,'' {\em IEEE Transactions on Automatic Control}, vol.~56, no.~9,
  pp.~2038--2049, 2011.

\bibitem{Bertsimas2004}
D.~Bertsimas, K.~Natarajan, and C.-P. Teo, ``Probabilistic combinatorial
  optimization: Moments, semidefinite programming, and asymptotic bounds,''
  {\em SIAM Journal on Optimization}, vol.~15, no.~1, pp.~185--209, 2004.

\bibitem{Bertsimas2006}
D.~Bertsimas, K.~Natarajan, and C.-P. Teo, ``Persistence in discrete
  optimization under data uncertainty,'' {\em Mathematical Programming},
  vol.~108, no.~2-3, pp.~251--274, 2006.

\bibitem{Natarajan2011}
K.~Natarajan, C.~P. Teo, and Z.~Zheng, ``Mixed 0-1 linear programs under
  objective uncertainty: A completely positive representation,'' {\em
  Operations Research}, vol.~59, no.~3, pp.~713--728, 2011.

\bibitem{Li2014}
X.~Li, K.~Natarajan, C.-P. Teo, and Z.~Zheng, ``Distributionally robust mixed
  integer linear programs: Persistency models with applications,'' {\em
  European Journal of Operational Research}, vol.~233, no.~3, pp.~459--473,
  2014.

\bibitem{Gavriel2012}
C.~Gavriel, G.~Hanasusanto, and D.~Kuhn, ``Risk-averse shortest path
  problems,'' in {\em 2012 IEEE 51st IEEE Conference on Decision and Control
  (CDC)}, pp.~2533--2538, IEEE, 2012.

\bibitem{Cheng2016}
J.~Cheng, J.~Leung, and A.~Lisser, ``New reformulations of distributionally
  robust shortest path problem,'' {\em Computers \& Operations Research},
  vol.~74, pp.~196--204, 2016.

\bibitem{Jaillet2016}
P.~Jaillet, J.~Qi, and M.~Sim, ``Routing optimization under uncertainty,'' {\em
  Operations Research}, vol.~64, no.~1, pp.~186--200, 2016.

\bibitem{Wang2020}
Z.~Wang, K.~You, S.~Song, and Y.~Zhang, ``Wasserstein distributionally robust
  shortest path problem,'' {\em European Journal of Operational Research},
  2020.

\bibitem{Zhu2018}
N.~Zhu, C.~Fu, and S.~Ma, ``Data-driven distributionally robust optimization
  approach for reliable travel-time-information-gain-oriented traffic sensor
  location model,'' {\em Transportation Research Part B: Methodological},
  vol.~113, pp.~91--120, 2018.

\bibitem{Dani2008}
V.~Dani, T.~P. Hayes, and S.~M. Kakade, ``Stochastic linear optimization under
  bandit feedback,'' in {\em 21st Annual Conference on Learning Theory},
  pp.~355--366, 2008.

\bibitem{Bubeck2012}
S.~Bubeck, N.~Cesa-Bianchi, and S.~M. Kakade, ``Towards minimax policies for
  online linear optimization with bandit feedback,'' in {\em Conference on
  Learning Theory} (S.~Mannor, N.~Srebro, and R.~C. Williamson, eds.),
  pp.~1--41, 2012.

\bibitem{Clare2012}
G.~Clare and A.~Richards, ``Air traffic flow management under uncertainty:
  application of chance constraints,'' in {\em Proceedings of the 2nd
  international conference on application and theory of automation in command
  and control systems}, pp.~20--26, IRIT Press, 2012.

\bibitem{Chen2017}
J.~Chen, L.~Chen, and D.~Sun, ``Air traffic flow management under uncertainty
  using chance-constrained optimization,'' {\em Transportation Research Part B:
  Methodological}, vol.~102, pp.~124--141, 2017.

\bibitem{Miller2000}
E.~D. Miller-Hooks and H.~S. Mahmassani, ``Least expected time paths in
  stochastic, time-varying transportation networks,'' {\em Transportation
  Science}, vol.~34, no.~2, pp.~198--215, 2000.

\bibitem{Fu1998}
L.~Fu and L.~R. Rilett, ``Expected shortest paths in dynamic and stochastic
  traffic networks,'' {\em Transportation Research Part B: Methodological},
  vol.~32, no.~7, pp.~499--516, 1998.

\bibitem{Hanasusanto2016}
G.~A. Hanasusanto, D.~Kuhn, and W.~Wiesemann, ``K-adaptability in two-stage
  distributionally robust binary programming,'' {\em Operations Research
  Letters}, vol.~44, no.~1, pp.~6--11, 2016.

\bibitem{Bickel2008}
P.~J. Bickel, E.~Levina, {\em et~al.}, ``Regularized estimation of large
  covariance matrices,'' {\em The Annals of Statistics}, vol.~36, no.~1,
  pp.~199--227, 2008.

\bibitem{Shapiro2001}
A.~Shapiro, ``On duality theory of conic linear problems,'' in {\em
  Semi-infinite programming}, pp.~135--165, Springer, 2001.

\bibitem{Boucheron2003}
S.~Boucheron, G.~Lugosi, and O.~Bousquet, ``Concentration inequalities,'' in
  {\em Summer School on Machine Learning}, pp.~208--240, Springer, 2003.

\bibitem{Sun2007}
J.~Sun, {\em The statistical analysis of interval-censored failure time data}.
\newblock Springer Science \& Business Media, 2007.

\bibitem{Ferson2007}
S.~Ferson, V.~Kreinovich, J.~Hajagos, W.~Oberkampf, and L.~Ginzburg,
  ``Experimental uncertainty estimation and statistics for data having interval
  uncertainty,'' {\em Sandia National Laboratories, Report SAND2007-0939},
  vol.~162, 2007.

\bibitem{Birge2011}
J.~R. Birge and F.~Louveaux, {\em Introduction to stochastic programming}.
\newblock Springer Science \& Business Media, 2011.

\bibitem{Berg1997}
M.~De~Berg, M.~Van~Kreveld, M.~Overmars, and O.~Schwarzkopf, ``Computational
  geometry,'' in {\em Computational geometry}, pp.~1--17, Springer, 1997.

\bibitem{Conforti2014}
M.~Conforti, G.~Cornu{\'e}jols, and G.~Zambelli, {\em Integer programming},
  vol.~271.
\newblock Springer, 2014.

\bibitem{Buchheim2018}
C.~Buchheim and J.~Kurtz, ``Robust combinatorial optimization under convex and
  discrete cost uncertainty,'' {\em EURO Journal on Computational
  Optimization}, vol.~6, no.~3, pp.~211--238, 2018.

\bibitem{Audet1997}
C.~Audet, P.~Hansen, B.~Jaumard, and G.~Savard, ``Links between linear bilevel
  and mixed 0--1 programming problems,'' {\em Journal of Optimization Theory
  and Applications}, vol.~93, no.~2, pp.~273--300, 1997.

\bibitem{Gupta2004}
A.~K. Gupta and S.~Nadarajah, {\em Handbook of beta distribution and its
  applications}.
\newblock CRC press, 2004.

\bibitem{Rockafellar2000}
R.~T. Rockafellar, S.~Uryasev, {\em et~al.}, ``Optimization of conditional
  value-at-risk,'' {\em Journal of risk}, vol.~2, pp.~21--42, 2000.

\bibitem{Goerigk2020}
M.~Goerigk, A.~Kasperski, and P.~Zieli{\'n}ski, ``Two-stage combinatorial
  optimization problems under risk,'' {\em Theoretical Computer Science},
  vol.~804, pp.~29--45, 2020.

\end{thebibliography}

\appendix
\section{Modeling the ambiguity set for moment-based formulation (\ref{DRO formulation 2})} \label{sec: app}
Here, we discuss how to construct the ambiguity set (\ref{ambiguity set 2}) from interval-censored observations $\xi^{(a)}_{k}$, $k \in \{1, \ldots, n_0\}$, $a \in A$. Specifically, for each baseline subinterval $[l^{(i)}_a, u^{(i)}_a]$, $i \in \mathcal{D}_a$, $a \in A$, the decision-maker learns whether a given observation $\xi^{(a)}_{k}$ for some $k \in \{1, \ldots, n_0\}$ belongs to this subinterval or not.

Henceforth, we fix $a \in A$ and derive an upper bound for the mean, $\widehat{\mu}_a$, for each $a \in A$. Estimation of $\widehat{\sigma}^2_a$ is analogous. Initially, we replace $\xi^{(a)}_{k}$ with its maximal possible value, $\overline{\xi}^{(a)}_{k}$, and assume that:
\begin{equation} \label{eq: estimate of mu}
\widehat{\mu}_a \approx \frac{1}{n_0}\sum_{k = 1}^{n_0} \overline{\xi}^{(a)}_{k}
\end{equation}

In particular, $\overline{\xi}^{(a)}_k$, $k \in \{1, \ldots, n_0\}$, can be defined as follows. Let
$$u^{(i^*)}_a := \min\{u^{(i)}_a: \xi^{(a)}_{k} \in [l^{(i)}_a, u^{(i)}_a], i \in \mathcal{D}_a\}$$
Then
\begin{equation} \label{eq: upper estimate for observations}
\overline{\xi}^{(a)}_k = \begin{cases}
\min\{l^{(i)}_a: \xi^{(a)}_{k} \leq l^{(i)}_a \mbox{ and } u^{(i^*)}_a \in [l^{(i)}_a, u^{(i)}_a], i \in \mathcal{D}_a\}, \mbox{ if } \mbox{the minimum exists} \\
u^{(i^*)}_a, \mbox{otherwise}
\end{cases}
\end{equation}

Similarly to the Hoeffding bounds in (\ref{eq: azuma-hoeffding 1}) and (\ref{eq: azuma-hoeffding 2}) we attempt to provide a probabilistic guarantee with respect to our estimate (\ref{eq: estimate of mu}). Specifically, we partition the support interval $[l_a, u_a]$ into a set of elementary subintervals $[L^{(j)}_a, U^{(j)}_a]$, $j \in \mathcal{W}_a$; see Section \ref{sec: DR SP without EC}. For each elementary subinterval we define Bernoulli random variables as follows:
\begin{equation} \nonumber
\chi_{kj}^{(a)} = \begin{cases}
1, \mbox{ if } \xi_k^{(a)} \in [L^{(j)}_a, U^{(j)}_a] \\
0, \mbox{ otherwise}
\end{cases}
\end{equation}
Furthermore, if $\chi_{kj}^{(a)} = 1$, then the maximal possible value of $\xi_k^{(a)}$ is uniquely defined from (\ref{eq: upper estimate for observations}). Thus, we define $\Xi_{j}^{(a)} := \overline{\xi}^{(a)}_k$, $j \in \mathcal{W}_a$, for $k$ such that $\chi_{kj}^{(a)} = 1$.

Then estimate (\ref{eq: estimate of mu}) can be expressed as:
\begin{equation} \nonumber
\widehat{\mu}_a \approx \frac{1}{n_0} \sum_{j \in \mathcal{W}_a} \Xi_{j}^{(a)} \sum_{k = 1}^{n_0} \chi_{kj}^{(a)}
\end{equation}
Observe that:
\begin{equation} \nonumber
\frac{1}{n_0} \sum_{j \in \mathcal{W}_a} \Xi_{j}^{(a)} \sum_{k = 1}^{n_0} \chi_{kj}^{(a)} - \mathbb{E} \{\frac{1}{n_0} \sum_{j \in \mathcal{W}_a} \Xi_{j}^{(a)} \sum_{k = 1}^{n_0} \chi_{kj}^{(a)}\} = \sum_{j \in \mathcal{W}_a}\Xi_{j}^{(a)} (\frac{1}{n_0}\sum_{k = 1}^{n_0} \chi_{kj}^{(a)} - q^{(a)}_j),
\end{equation}
where $q^{(a)}_j = \mathbb{Q}^0_a\{c_a \in [L^{(j)}_a, U^{(j)}_a]\}$, $j \in \mathcal{W}_a$.

Hoeffding inequality implies that for any given $\varepsilon_a > 0$ the following inequality holds:
\begin{equation} \nonumber
q^{(a)}_j - \frac{1}{n_0}\sum_{k = 1}^{n_0} \chi_{kj}^{(a)} \geq \varepsilon'_{aj} := \frac{\varepsilon_a}{W_a \Xi_{j}^{(a)}},
\end{equation}
with probability of at most $\exp(-2r{\varepsilon'}_{aj}^{2})$. By leveraging Bonferroni's inequality we have:
\begin{equation} \label{eq: confidence interval for mu}
\widehat{\mu}_a = \sum_{j \in \mathcal{W}_a}\Xi_{j}^{(a)} q^{(a)}_j \leq \sum_{j \in \mathcal{W}_a}\Xi_{j}^{(a)} \frac{1}{n_0}\sum_{k = 1}^{n_0} \chi_{kj}^{(a)} + \varepsilon_a = \frac{1}{n_0}\sum_{k = 1}^{n_0} \overline{\xi}^{(a)}_{k} + \varepsilon_a,
\end{equation}
with probability of at least $1 - \sum_{j \in \mathcal{W}_a} \exp(-2r{\varepsilon'}_{aj}^{2})$. Finally, by setting
$$\eta_0 = 2|A| \sum_{j \in \mathcal{W}_a} \exp(-2r{\varepsilon'}_{aj}^{2}) = 2|A| \sum_{j \in \mathcal{W}_a} \exp(-2r \Big(\frac{\varepsilon_a}{W_a \Xi_{j}^{(a)}}\Big)^{2})$$ for each particular $a \in A$ we find the value of $\varepsilon_a$, which provides a prescribed confidence level, $1 - \eta_0$, for the ambiguity set (\ref{ambiguity set 2}).

Despite the fact that the first-order moment constraint in (\ref{ambiguity set 2}) is unique for each $a \in A$ its probability of violation is overestimated, i.e., violation of (\ref{eq: confidence interval for mu}) does not necessarily imply violation of the constraint $\mathbb{E}_{\mathbb{Q}_a}\{c_a\} \leq \widehat{\mu}_a$. In fact, concentration inequalities cannot be applied directly to the moment constraints whenever the data is interval-censored. Alternatively, if random samples $\xi_k^{(a)}$, $k \in \{1, \ldots, n_0\}$, are available to the decision-maker for each particular arc $a \in A$ (and $n_0$ exceeds some instance-dependent threshold), then a confidence region for the mean and covariance matrix can be constructed using techniques from \cite{Delage2010}.

Overall, the outlined approach is myopic and not without limitations. Nevertheless, we attempt to exploit the key properties of our formulation (\ref{DRO formulation}) to provide a valid comparison of two distributionally robust models, namely, (\ref{DRO formulation}) and (\ref{DRO formulation 2}).

\end{document}